\documentclass[12pt,reqno]{amsart}
\RequirePackage{multirow}
\usepackage{xifthen,setspace,multirow}
\usepackage{diagbox,graphicx}
\usepackage{caption,subcaption}
\usepackage{float}
\usepackage[margin=3cm]{geometry}
\usepackage{bbm}
\usepackage{xcolor}
\usepackage{enumitem}
\usepackage{bm}
\usepackage[linesnumbered,ruled,vlined]{algorithm2e}
\usepackage{hyperref}
\usepackage{esvect}
\usepackage{amsmath,amssymb,mathrsfs} % provides numberwithin (and lots more)
\usepackage{lipsum}  % for sample text

\newtheorem{remark}{Remark}
\newtheorem{theorem}{Theorem}
\newtheorem{corollary}{Corollary}

\newtheorem{lem}{Lemma}
\newtheorem{assumption}{Assumption}

\newcommand{\vertiii}[1]{{\left\vert\kern-0.25ex\left\vert\kern-0.25ex\left\vert #1 
    \right\vert\kern-0.25ex\right\vert\kern-0.25ex\right\vert}}

\renewcommand{\leq}{\leqslant}
\renewcommand{\geq}{\geqslant}

\newcommand{\bpk}{{\binom{p-1}{k-1}}}

\newcommand{\p}{{\mathbb{P}}}

\allowdisplaybreaks

\title[Tensor Recovery in High-Dimensional Ising Models]{Tensor Recovery in High-Dimensional Ising Models}

\author[Liu]{Tianyu Liu}
\address{email {\tt tianyu.liu@u.nus.edu}}

\author[Mukherjee]{Somabha Mukherjee}
\address{email {\tt somabha@nus.edu.sg}}

\author[Biswas]{Rahul Biswas}
\address{email {\tt rbiswas1@uw.edu}}

\begin{document}

\begin{abstract} 
The $k$-tensor Ising model is an exponential family on a $p$-dimensional binary hypercube for modeling dependent binary data, where the sufficient statistic consists of all $k$-fold products of the observations, and the parameter is an unknown $k$-fold tensor, designed to capture higher-order interactions between the binary variables. In this paper, we describe an approach based on a penalization technique that helps us recover the signed support of the tensor parameter with high probability, assuming that no entry of the true tensor is too close to zero. The method is based on an $\ell_1$-regularized node-wise logistic regression, that recovers the signed neighborhood of each node with high probability. Our analysis is carried out in the high-dimensional regime, that allows the dimension $p$ of the Ising model, as well as the interaction factor $k$ to potentially grow to $\infty$ with the sample size $n$. We show that if the minimum interaction strength is not too small, then consistent recovery of the entire signed support is possible if one takes $n = \Omega((k!)^8 d^3 \log \bpk)$ samples, where $d$ denotes the maximum degree of the hypernetwork in question. Our results are validated in two simulation settings, and applied on a real neurobiological dataset consisting of multi-array electro-physiological recordings from the mouse visual cortex, to model higher-order interactions between the brain regions.
\end{abstract}

	%\spacing{1.125}


	\maketitle

\section{Introduction}\label{int}
The Ising model, a discrete exponential family for modeling dependent binary data, was initially used by physicists as a model for ferromagnetism \cite{ising}. Since then, this model was applied immensely in diverse fields such as computational biology, neural networks, social sciences, image processing, spatial statistics and election forecasting \cite{spatial,cd_trees,geman_graffinge,disease,neural,innovations,elhj1,isinggenapl5}. The classical $2$-spin Ising model is an exponential family on the binary hypercube, whose sufficient statistic involves all pairwise products of the binary ($\pm 1$-valued) observations, and whose parameter is an interaction matrix, designed to capture pairwise interactions between the binary variables. The problem of structure learning in Ising models relates to estimating the interaction matrix, given access to multiple i.i.d. samples from the same model. A significant amount of work has been done in the literature on structure recovery in classical $2$-spin Ising models, the notable ones being \cite{structure_learning,bresler,discrete_tree,graphical_models_algorithmic,wainwright,graphical_models_binary,lokhov}. Daskalakis et al. \cite{cd_testing} studied the problems of identity and independence testing, and Neykov et al. \cite{high_tempferro,neykovliu_property} considered, instead of recovering the full structure, the more fundamental problem of graph property testing, such as connectivity, presence of cycles and maximum clique size, given access to multiple samples from an Ising model. 

A different branch of research on estimation in Ising models, assumes that the interation matrix is known upto some scale factor, and sometimes assumes the presence of external magnetic fields in the model. The main focus in this area is to estimate the scale factor of the interaction matrix (referred to as \textit{inverse temperature} in statistical physics), and the external magnetic fields \cite{chatterjee,pg_sm,BM16}, and to prove asymptotics of these estimators for some standard Ising models \cite{comets}. The techniques used in this field are quite different from the ones used in the literature of structure learning, primarily due to the fact that estimation in the former area is often based on only one sample, in contrast with the multiple-sample regime considered in the latter.

Structure learning in Ising models has applications in diverse disciplines, such as epidemic network modeling, statistical physics, image processing, machine learning and spatial transcriptomics. For example, in a contagious epidemic network, of utmost importance is the network effect (probably more than the personal attributes such as age, weight, immunity, smoking habits, etc.), i.e. other people in the network that a particular person came into contact with. Another area where structure learning is highly relevant, is the field of spatial transcriptomics. This is a relatively modern area in biology, that uses technologies designed for vastly parallelized measurement of cell transcriptomes in situ. In contrast to single cell sequencing, spatial transcriptomics retains information regarding the spatial arrangement of the cells, which can be thought to be encoded into a Voronoi neighborhood graph, with the nodes denoting the cells and edges being drawn between proximally located cells. For each node, the genetic/protein expressions are recorded, and the typical goal is to understand how the spatial structure of the cells contribute to their phenotypes. In neuroscience, learning interaction structures between neurons is a popular subject \cite{biscfc1}, and Ising models have been used to model and infer interactions between spiking activity of a population of neurons \cite{isingfc_1, isingfc_2}, and collective properties of the neuronal network \cite{isingneuro_1, isingneuro_2}.%, and  \textcolor{magenta}{Few lines about the neuron setting/application.}

However, in most real-world scenarios, pairwise interactions are not enough to capture the complex dependencies arising in a network structure, but one has to take into account higher order peer-group effects. To elaborate, it is often more reasonable for an individual to choose a binary attribute if many of his/her friends have also chosen the same. Another example comes from chemistry, where it is known that the atoms on a crystal surface (adatoms) do not interact just in pairs, but in triangles, quadruplates and higher order tuples. In neuroscience, it is well known that a single neuron receives and sends impulses from and to multiple neurons, thereby motivating the need to consider higher order interactions. Hypergraphs/ interaction tensors are useful ways of quantifying higher-order relational data arising naturally in a wide variety of applications \cite{hypergraph_learning,hypergraph_applications,hypergraph_complex_data,hypergraph_image,hypergraph_multimedia,hypergraph_gene}, and in order to understand the complex relationships of the variables in such datasets, one natural choice is to consider tensor Ising models \cite{barra, oliveira, gardner, bovier, mukherjeeestimation_b}, where the interaction matrix is replaced by a tensor, encoding the strength of the interactions between, not just pairwise, but groups or tuples of individuals. Estimating the support of this unknown tensor is of natural interest, and in this paper, we achieve this by running penalized node-wise logistic regressions that recover the signed neighborhoods of each node with high probability. The point to note is that one cannnot use the model likelihood function here to do this, because the corresponding normalizing constant is inexplicit and computationally intractable. A computationally efficient alternative is to work with the pseudolikelihood function \cite{chatterjee, besag_lattice, besag_nl}, which for every node, computes the product of the conditional distributions of the observation at that node given all the remaining nodes, over all the samples available. Unlike the likelihood function, the pseudolikelihood is free of the intractable normalizing constant, and is in fact, computationally explicit.

\subsection{The Tensor Recovery Problem} 
The $k$-tensor Ising model (see \cite{barra, oliveira, gardner, bovier, mukherjeeestimation_b}) is a probability distribution on the set $\{-1,1\}^p$, defined as:
\begin{equation}\label{modeldef}
    \p_{\bm J}(\bm x) := \frac{1}{Z(\bm J)}e^{H(\bm x)}\quad(\bm x \in \{-1,1\}^p)
\end{equation}
where $\bm J := ((J_{r_1,\ldots,r_k}))_{(r_1,\ldots,r_k)\in [p]^k}$ denotes a $k$-fold tensor with $[p] := \{1,\ldots,p\}$ and $$H(\bm x) := \sum_{(r_1,\ldots,r_k)\in [p]^k} J_{r_1,\ldots,r_k} x_{r_1}\ldots x_{r_k}.$$
Hereafter, we will assume that the tensor $\bm J$ satisfies the following properties:
\begin{enumerate}
    \item $\bm J$ is symmetric, i.e., $J_{r_1,...,r_k}= J_{r_{\sigma (1)},...,r_{\sigma(k)}}$ for every $(r_1,\ldots,r_k)\in [p]^k$ and every permutation $\sigma$ of $\{1,...,k\}$,
    \item $\bm J$ has zeros on the \textit{diagonals}, i.e., $\bm J_{r_1...r_k}=0$, if $r_s=r_t$ for some $1\leqslant s < t \leqslant k$.
\end{enumerate}

We will also assume that $p \ge 4$ and $2\le k\le p-1$ for technical reasons that will become clear later. Suppose that we are given a collection $\mathfrak{X}^{n}:=$ $\left\{\bm x^{(1)}, \ldots, \bm x^{(n)}\right\}$ of $n$ samples from the model \eqref{modeldef}. Our aim is to infer the underlying tensor $\bm J$ based on this sample $\mathfrak{X}^n$. A common example of such a tensor $\bm J$ is the adjacency of a $k$-uniform hypergraph. Analogous to the principal goal of \textit{graphical model selection}, a natural aim in this setup is to recover the hyperedge set of the tensor $\bm J$. In this article, we consider the slightly stronger problem of \textit{signed hyperedge recovery}. To elaborate, we define the signed-edge tensor corresponding to $\bm J$ as:
$$\bm J^* := ((\mathrm{sgn}(J_{r_1,\ldots,r_k})))_{(r_1,\ldots,r_k)\in [p]^k}$$ where $\mathrm{sgn}(t) := t/|t|$ (if $t\ne 0$) and $\mathrm{sgn}(0) := 0$.
Following the idea in \cite{wainwright}, we apply an $\ell_1$-penalized LASSO approach to recover the signed-edge tensor $\bm J^*$. 
Our theoretical results focus on showing consistency of the signed hyperedge recovery algorithm, for which we go beyond the classical statistical framework of fixed $p$ and $n\rightarrow \infty$, and work under a high-dimensional setting, where both $p$ and $k$ are allowed to grow with $n$. Moreover, if we define:
$$d_r := \sum_{1\le r_1<\ldots<r_{k-1}\le p} \mathbbm{1}_{\{J_{r,r_1,\ldots,r_{k-1}} \ne 0\}}\quad\text{and}\quad d := \max_{r\in [p]} d_r~,$$ then we also allow the maximum degree $d$ to grow with $n$ in our framework. The precise relations between these four quantities $n,p,k$ and $d$ in order to guarantee consistent recovery of $\bm J^*$ are specified in the theoretical results in Section \ref{thres}.

Recovering the signed-edge tensor $\bm J^*$ is equivalent to recovering for each vertex $r$, the vector:
$$\bm J_r^* := ((J_{r,r_1,\ldots,r_{k-1}}^*))_{(r_1,\ldots,r_{k-1})\in T_r}$$ where $T_r := \{(r_1,\ldots,r_{k-1}) \in ([p]\setminus \{r\})^{k-1}: 1\le r_1<\ldots<r_{k-1} \le p\}$. For this, we implement the following node-wise $\ell_1$-regularized pseudolikelihood approach:
\begin{equation}\label{lasso9}
    \min_{\bm J_r \in \mathbb{R}^{\binom{p-1}{k-1}}} \ell(\bm J_r;\mathfrak{X}^n) + \lambda \|\bm J_r\|_1
\end{equation}
where 
$$\ell(\bm J_r;\mathfrak{X}^n) := -\frac{1}{n} \sum_{i=1}^n \log \p_{\bm J}(x_r^{(i)}| \bm x_{\setminus r}^{(i)})$$
and $\bm x_{\setminus r}^{(i)} := (x_t^{(i)})_{t\ne r}$. A straightforward computation shows that:
$$\p_{\bm J}(x_r| \bm x_{\setminus r}) = \frac{\exp(k x_r m_r(\bm x))}{2\cosh(k x_r m_r(\bm x))}$$
where $m_r(\bm x) := \sum_{(r_1,\ldots,r_{k-1}) \in [p]^{k-1}} J_{r,r_1,\ldots,r_{k-1}} x_{r_1}\ldots x_{r_{k-1}}$.

We define the \textit{hyperedge set} as $E := \{\bm e \subseteq [p]: J_{\bm e} \ne 0\}$ and the \textit{neighborhood} and \textit{signed neighborhood} of each vertex $r$ respectively as:
$$\mathcal{N}(r) := \{\bm e' \subseteq [p]: |\bm e'| = k-1, \{r\}\cup \bm e' \in E\}\quad\text{and}\quad \mathcal{N}^{\pm}(r) := \{J_{\{r\}\cup \bm e'}^* ~\bm e': \bm e' \in \mathcal{N}(r)\}~.$$
Note that $\mathcal{N}^\pm (r)$ can be recovered for every vertex $r$ if we can infer the vector $\bm J_r^*$, for which we solve the optimization program \eqref{lasso9}.

\subsection{Assumptions}
We require certain assumptions to ensure that our method works consistently. To state these assumptions, we need a few preliminary notations.
For any fixed node $r\in V$, we define a $\binom{p-1}{k-1}\times \binom{p-1}{k-1}$ matrix of the form:
\begin{equation}\label{qstar}
    \bm Q_r:=-\mathbb{E}_{\bm J}\left[\nabla_{\bm J_r}^2\log \mathbb{P}_{\bm J}(X_r|\bm X_{\backslash r})\right]
\end{equation}
which can be written more explicitly as:
\begin{equation*}
    \bm Q_r:=\mathbb{E}_{\bm J}\left[\eta_r(\bm X;J) \bm X_{\cdot r} \bm X_{\cdot r}^T\right],
\end{equation*}
where
\begin{equation*}
   \eta_r(\bm X;\bm J):=\frac{4(k!)^2 e^{2k X_r m_r(\bm X)}}{(e^{2k X_r m_r(\bm X)}+1)^2},
\end{equation*}
and the $\binom{p-1}{k-1}$ dimensional vector $\bm X_{\cdot r}$ is defined as 
\begin{equation*}
    \bm X_{\cdot r}:=(X_{r_1}\cdots X_{r_{k-1}})_{(r_1,\ldots,r_{k-1})\in T_r}~.
\end{equation*}
When the node $r$ is clear from the context and there is no scope of any confusion, we will henceforth abbreviate $\bm Q_r$ by $\bm Q$. Besides, we will denote the set of all hyperedges containing $r$ by $S_r$ (also abbreviated by $S$ when there is no scope of confusion), i.e.
$$S_r := \{\bm e\in E: r\in \bm e\}.$$
 Following this notation we define $\bm Q_{S_r S_r}$ as the $d_r\times d_r$ sub-matrix of $\bm Q$ indexed by $S_r$. With the above notations, we are now ready to state our assumptions.

 \begin{assumption}[Dependency Condition]\label{depcon}
There exist constants $C_{\min}$, $D_{\max}$ such that
\begin{equation}\label{a1}
    \begin{aligned}
    \Lambda_{\min}(\bm Q_{SS})&\geqslant C_{\min},\\
    \Lambda_{\max}(\mathbb{E}_{\bm J}[\bm X_{\cdot r} \bm X_{\cdot r}^T])&\leqslant D_{\max},
    \end{aligned}
\end{equation}
 \end{assumption}

  The first condition in \eqref{a1} bounds the minimum eigenvalue of the Fisher information matrix corresponding to the subset of relevant covariates, and the second condition in \eqref{a1} states that the relevant covariates are not overly dependent. 

  \begin{assumption}[Incoherence Condition]\label{a27}
   There exists an $\alpha \in (0,1]$ such that
\begin{equation}\label{a2}
    \vertiii{\bm Q_{S^c S}\bm Q_{SS}^{-1}}_{\infty} \leqslant 1-\alpha,
\end{equation} 
where $\vertiii{\cdot}_\infty$ refers to the matrix $\ell^\infty$ norm.
  \end{assumption}
  
  Condition \eqref{a2} restraints the influence of the irrelevant covariates on all the relevant covariates. 

\subsection{Organization} The rest of the paper is organized as follows. In Section \ref{thres}, we state the main theoretical results in this paper, on consistent tensor recovery, and give a brief sketch of the proof. Section \ref{numstudy} is devoted to applications of the recovery algorithm on some experimental and real-life neurobiological datasets. The simulation framework includes two different scenarios, one where samples are generated from Ising models on regular hypergraphs, and in another, where data is simulated from a 3-tensor Ising model on a user friendship network obtained from the Last.fm dataset. In the real-life data analysis section, the method is applied on a neurobiological dataset consisting of electro-physiological recordings from the visual cortex region in the mouse brain. The proofs of the main results are given in Section \ref{mainprff}. In Section \ref{discuss_sec}, we summarize our main contributions, and provide directions for future research. Proofs of some technical lemmas are given in the appendix.

\section{Theoretical Results}\label{thres}

In this section, we state the main theoretical results of this paper. Specifically, we give sufficient conditions on the tuple $(n,p,d,k)$ and the regularization parameter $\lambda_n$, that guarantee successful recovery of the signed neighborhood vectors.

\begin{theorem}\label{thm1}
Suppose that the regularization parameter $\lambda$ is chosen as:
\begin{equation}\label{reg}
    \lambda \propto k! \sqrt{\frac{\log \bpk}{n}} ~.
\end{equation}
Then there exists a positive constant $L$ independent of $(n, p, d, k)$, such that if
\begin{equation}\label{scaling}
   n >L (k!)^8 d^3 \log \binom{p-1}{k-1}
\end{equation}
then for each node $r$, the following properties hold with probability at least $1-M \exp\left(-h(n,d,k)\right)$ for some constant $M>0$, where $$h(n,d,k) := \min\left\{\log \bpk~,~ \frac{K_1n}{d^3(k!)^4} -K_2\log d -\log\left[\bpk - d\right]  \right\}$$
for some constants $K_1,K_2>0$.
\vspace{0.1in}

\noindent (a) The $\ell_1$-penalized logistic regression (\ref{lasso9}) has a unique solution, and hence uniquely specifies an estimated signed neighborhood $\widehat{\mathcal{N}}^\pm (r)$.
\vspace{0.1in}

\noindent (b) The estimated signed neighborhood $\widehat{\mathcal{N}}^\pm (r)$ correctly excludes all hyperedges not in the true neighborhood. Moreover, it correctly includes all hyperedges $\bm e$ containing $r$, for which $|J_{\bm e}| \geqslant \frac{10}{C_{\min }} \lambda \sqrt{d} $.
\end{theorem}

The following result is an easy consequence of Theorem \ref{thm1} by a further union bound applied on all the nodes of the hypergraph.

\begin{corollary}\label{cor117}
Suppose that in addition to Conditions \eqref{reg} and \eqref{scaling}, we have the following two conditions:
\vspace{0.1in}

\noindent (a) $p \rightarrow \infty$.
\vspace{0.1in}

\noindent (b) $\min_{\bm e \in E}|J_{\bm e}| \geqslant \frac{10}{C_{\min}} \lambda \sqrt{d}$ for sufficiently large $n$.
\vspace{0.1in}

 \noindent Then our model selection algorithm is consistent, i.e. if $\widehat{E}$ denotes the estimated hyperedge set, then
$$\p\left(\widehat{E} = E\right) \rightarrow 1~\text{as}~n\rightarrow \infty~.$$ 
\end{corollary}

\subsection{Sketch of Proof}
The proof adopts the techniques in \cite{wainwright}, modulo some modifications taking care of the tensor case. To begin with, note that the optimization problem \eqref{lasso9} can be re-written as:
\begin{equation}\label{min}
    \min_{\bm J_r \in \mathbb{R}^{\bpk}} \frac{1}{n} \sum_{i=1}^{n} f\left(\bm J_r ; \bm x^{(i)}\right)-k\sum_{(r_1,...,r_{k-1}) \in [p]\setminus \{r\}} J_{r, r_1,\ldots, r_{k-1}} \widehat{\mu}_{r, r_1, \ldots, r_{k-1}}+\lambda\|\bm J_r\|_{1}
\end{equation}
where
\begin{equation*}
    f(\bm J_r;\bm x):=\text{log}\left(e^{k m_r(\bm x)} +e^{-km_r(\bm x)} \right)
\end{equation*}
and $\widehat{\mu}_{r, r_1,\ldots,r_{k-1}}:= \frac{1}{n}\sum\limits_{i=1}^n x^{(i)}_r x^{(i)}_{r_1}\cdots x^{(i)}_{r_{k-1}}$ are the empirical moments. The proof uses a primal-dual witness approach. The first step is to show that any optimal primal solution to \eqref{min} must vanish on all index tuples which are not hyperedges, and under standard assumptions, a unique optimal solution $\hat{\bm J}_r$ is guaranteed. The next step is to provide an upper bound to the term $\|\nabla^2 \ell(\bm J_r,\mathfrak{X}^n)(\hat{\bm J}_r - \bm J_r)\|_\infty$, and a lower bound on the minimum eigenvalue of $\nabla^2\ell(\bm J_{SS},\mathfrak{X}^n)$. These together, will establish consistency of $\hat{\bm J}_r$ towards $\bm J$. Further details are provided in Section \ref{mainprff}.

\begin{remark}\label{lucrem}
Instead of learning the Ising tensor structure, we can instead ask the question that how closely can we learn the entire Ising distribution \eqref{modeldef}, in terms of a certain notion of distance between measures. It follows from the proofs of Theorem 1.1 and Theorem 1.3 in \cite {luc} that the learning rate of the tensor Ising model \eqref{modeldef} in terms of the total variation distance, is bounded above by a constant factor times $\min\{1,\sqrt{|E|/p}\}$, where $E$ denotes the edge set of the underlying hypernetwork. It is also proved in \cite{luc} that this rate is minimax optimal for $k=2$. The proof of the lower bound in \cite{luc} involves concentration inequalities for quadratic forms, and hence, is not expected to adapt as it is for the tensor Ising models. This can be a potential direction for future research.
\end{remark}

\section{Numerical Study}\label{numstudy}
This section is devoted to applying the tensor recovery algorithm proposed in Section \ref{int} to some simulated and real-life scenarios. 
\subsection{Simultation Study}\label{simstudy}
In this section, we present some numerical experiments that illustrate the performance of the tensor recovery algorithm. In these examples, we study the recovery rate of the algorithm based on samples simulated from Ising models on sparse hypergraphs. The built-in lib-linear solver in python ``\textit{sklearn}" package was used to solve the $\ell_1$-regularized logistic regression in all cases. 
\subsubsection{Regular hypergraphs}
We analyzed the performance of our algorithm on samples generated from Ising models on regular hypergraphs. The coefficients of the hyperedges are assigned sign $+1$, and samples are then simulated from the corresponding tensor Ising model by a Gibbs sampler, with the aim of inferring the hyperedge signs. The recovery rate was then examined under the following scaling:
$$n = \frac{\alpha}{6\times 10^6} (k!)^8 d^3 \log \binom{p-1}{k-1}$$ where as usual, $n$ denotes the sample size and $p$ denotes the number of nodes in the $k$-uniform, $d$-regular hypergraph. The scaling parameter $\alpha$ ranged from $0.2$ to $2$, with graph sizes $p\in\{32, 64, 128\}$, $d$ is set at $3$, and the cardinality of hyperedges is fixed at $k=3$. The regularization parameter $\lambda$ was set as $\lambda = c \sqrt{k\log p/n}$, where $c$ is a constant factor tuned according to the Bayesian Information Criterion. For each node $r$, the optimal value of $\lambda$ is tuned by minimizing the BIC value with $LassoLarsIC$ method in the $``sklearn"$ package. Then the average of all $\lambda$s is taken as the regularization parameter. 

\begin{figure}[h]
    \centering
    \includegraphics[width=0.5\textwidth]{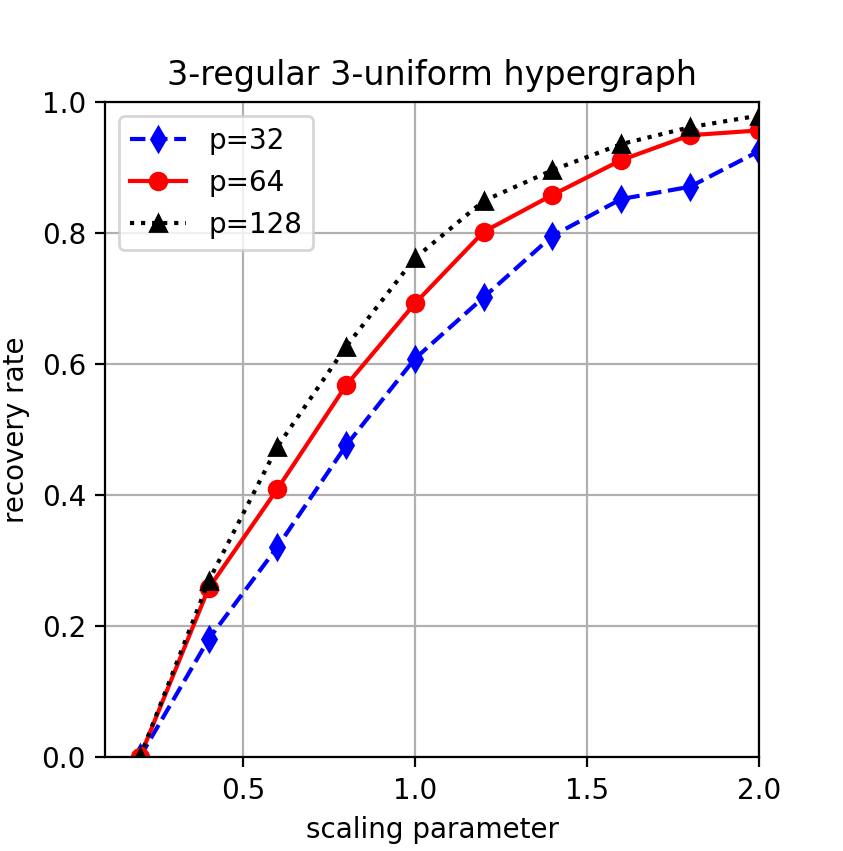}
    \caption{Plots of the recovery rate against the scaling parameter $\alpha = (6\times 10^6) n/[(k!)^8 d^3 \log \binom{p-1}{k-1}]$}
    \label{fig1}
\end{figure}

\begin{figure}[h]
    \centering
    \includegraphics[width=0.5\textwidth]{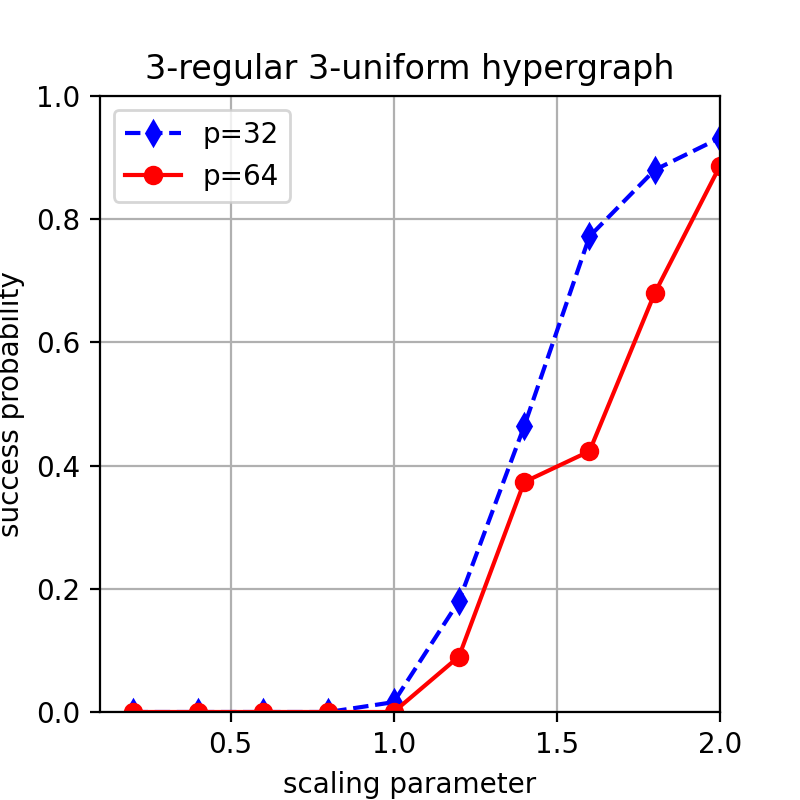}
    \caption{Plots of the probability of success against the scaling parameter $\alpha = (1.5 \times 10^6) n/[(k!)^8 d^3 \log \binom{p-1}{k-1}]$}
    \label{fig2}
\end{figure}

Figure \ref{fig1} shows the recovery rate versus the scaling parameter as a factor of the sample size, with $k=3$. Each curve corresponds to a given graph size $p\in\{32, 64, 128\}$. The recovery rate is defined as:
$$r = \frac{|\{\bm e \in E: \widehat{J}_{\bm e}^* = J_{\bm e}^*\}|}{|E|}~.$$ For each graph size and each scaling parameter, $50$ trials were conducted and an average recovery rate was computed. It can be seen that the three curves have similar shapes, despite the difference in graph size. Moreover, with increase in the scaling parameter (or equivalently, with increase in the sample size), the recovery rate approaches $1$, which illustrates the validity of Theorem \ref{thm1} and Corollary \ref{cor117}.

In Figure \ref{fig2}, we plot the success probability of complete recovery for graph sizes $32$ and $64$ against the scaling parameter $\alpha$ as a function of the sample size, where the success probability of complete recovery is defined as the fraction of cases where the edge set is completely recovered in a number of repeated independent trials of the algorithm. To be specific, we take:
$$\alpha = \frac{(1.5 \times 10^6) n}{(k!)^8 d^3 \log \binom{p-1}{k-1}}$$
The success probability also approaches $1$ as the scaling parameter $\alpha$ increases.

\subsubsection{The Last.fm Dataset}

\begin{figure}[h]
    \centering
    \includegraphics[width=0.5\textwidth]{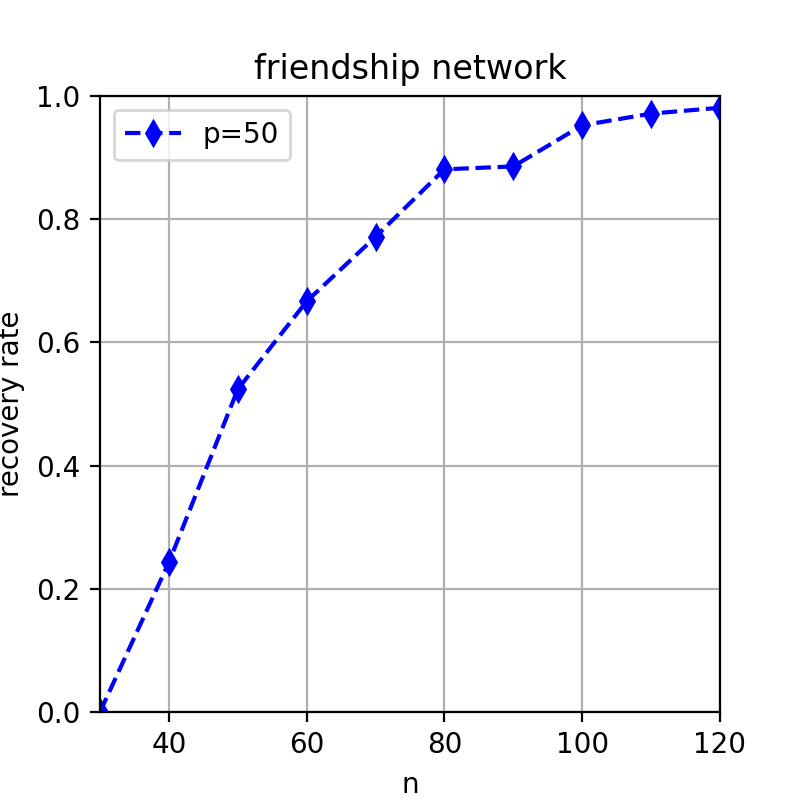}
    \caption{Plot of the recovery rate against the number of samples n, for $3$-uniform hypergraph truncated from the friendship network, with $k=3$}
    \label{fig3}
\end{figure}

\begin{figure}[h]
    \centering
    \includegraphics[width=0.5\textwidth]{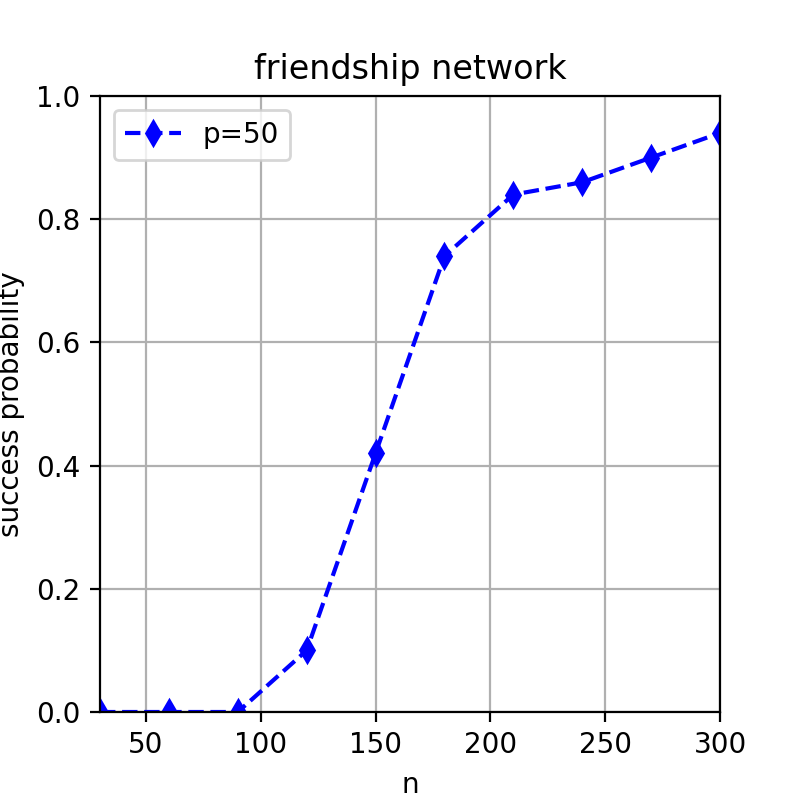}
    \caption{Plot of the success probability against the number of samples n, for $3$-uniform hypergraph truncated from the friendship network, with $k=3$}
    \label{fig4}
\end{figure}

The Last.fm dataset (\url{http://millionsongdataset.com/lastfm/}), which is a part of the Million Song Dataset (\url{http://millionsongdataset.com/}) consists of a list of 1892 users, their friendship network, and their most favorite artists (see \cite{fiftyfour,seventythree}). It was shown in \cite{dissertation} that users' preference for music artists are influenced by possible peer group effects present in the social network between the users, and hence, a tensor Ising model might be a good fit. In this section, we extracted the triangles from the user friendship network to create a $3$-uniform tensor, and fitted an Ising model on this tensor. We generated $n \in [30,300]$ samples from this tensor Ising model, based on which, the hypergraph structure was inferred. Figure \ref{fig3} shows the recovery rate as a function of the number of samples. Figure \ref{fig4} shows the success probability as a function of the number of samples. Once again, this approaches $1$ with increase in the sample size.

\subsection{Applications in a Neurobiological Dataset}\label{realstudy}
\begin{figure}[t!]
    \centering
    \includegraphics[width = \textwidth]{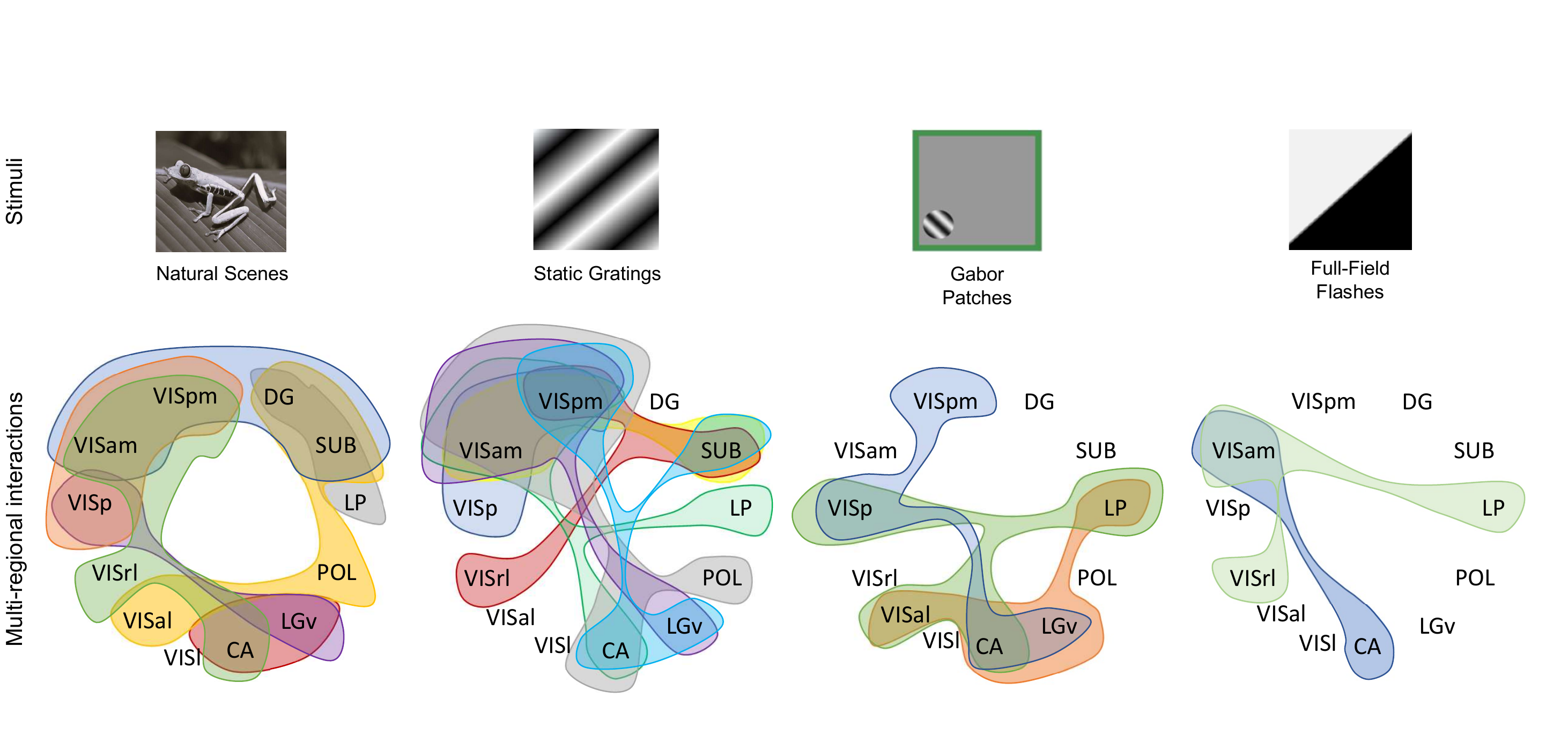}
    \caption{Estimated multi-body interactions between different regions of the mouse brain.}
    \label{fig:4body}
\end{figure}

We implement our methods on neural data acquired through electrophysiological recordings from the Visual Coding Neuropixels dataset of the Allen Brain Observatory \cite{de2020large}. We limit our study to a male mouse aged 116 days (Session ID 791319847) with 555 neurons whose spike trains were recorded simultaneously via six Neuropixel probes. The spike trains were recorded at a frequency of 1 KHz throughout the entirety of the experiment. Our analysis focuses on the following four stimulus categories (see Figure \ref{fig:4body}) \cite{biscfc2}.

\begin{enumerate}
    \item     Natural scenes, consisting of 118 images from three databases presented briefly for 250ms each, repeated 50 times in random order with intermittent blank intervals.
    \item     Static gratings, consisting of full-field sinusoidal gratings with varying orientations, spatial frequencies, and phases, resulting in 120 conditions. Each condition is presented briefly (250 ms) before being replaced with a different condition, repeated 50 times in random order with intermittent blank intervals.
    \item     Gabor patches with 3 orientations presented at different points in a 9 $\times$ 9 visual field, each presented for 250ms and repeated 50 times in random order with intermittent blank intervals.
    \item     Full-field flashes, lasting for 250ms followed by a blank interval of 1.75s and repeated 150 times.
\end{enumerate}

These stimuli range from natural scenes to artificial stimuli, and the purpose of this study is to explore how different stimuli patterns affect multi-neuron interactions. Dynamic stimuli such as natural movies and drifting gratings are excluded from the analysis as they require further investigation and interpretation in future studies.

To preprocess the data, we converted the spike trains, which were recorded at a frequency of 1 KHz, to a bin size of 10 ms. We achieved this by aggregating and separating them based on the start and end times of each stimulus presentation, resulting in Peri-Stimulus Time Histograms (PSTH) with a bin size of 10 ms. To obtain a smoothed version of the PSTH for each neuron and each stimulus presentation, we used a Gaussian smoothing kernel with a bandwidth of 16ms. For each stimulus presentation, we used the smoothed PSTHs as input to infer the functional connectivity (FC) between the neurons. We selected the set of neurons that were active in at least 25\% of the bins in the PSTH for each stimulus presentation and collected the set of unique neurons across all stimuli. This resulted in 33, 27, 12 and 23 active neurons for natural scenes, static gratings, Gabor patches, and flashes, respectively, and a total of 44 unique active neurons overall. We further estimated the numerical derivatives of the PSTH for each active neuron and obtained a sequence with $+1$, if the derivative is positive, and $-1$, if the derivative is negative. We then abandoned the points where that derivative is $0$ and obtain 3784 out of 5950 recordings for natural scenes, 4226 out of 6000 recordings for static gratings, 2945 out of 3645 recordings for gabors and 94 out of 150 recordings for flashes. To reduce the dependence from the time series, we picked 1 in every 3 time points. For flashes, we took all the time points as samples.

%We recover the $4$-neuron interactions using a $4$-fold tensor Ising model, and the neuronal interactions are recovered for each stimulus type and compared across the different stimulus types.

We recover the $4$-neuron interactions using a $4$-tensor Ising model, and the neuronal interactions are recovered for each stimulus type and compared across the different stimulus types. Figure \ref{fig:4body} shows examples of such estimated multi-body interactions, after labeling the neurons based on their brain regions. For example, in the natural scenes scenario, some multi-body interactions include VISam-VISpm-DG-SUB, VISp-VISam-VISpm, VISam-VISpm-VISrl-CA, and VISal-POL-SUB-DG. In static gratings, the typical multi-body interactions noticed, are VISam-VISpm-CA-POL, VISpm-SUB-VISrl, VISam-VISpm-LGv, and VISam-VISpm-CA-POL. In Gabors, VISpm-VISp-CA-LGv, VISp-VISal-CA-LP, and VISal-CA-LGv-LP comprise some of the interactions, while in Flashes, some of these interactions include VISam-VISrl-LP and VISam-CA.

\begin{figure}[t!]
    \centering
    \includegraphics[width = \textwidth]{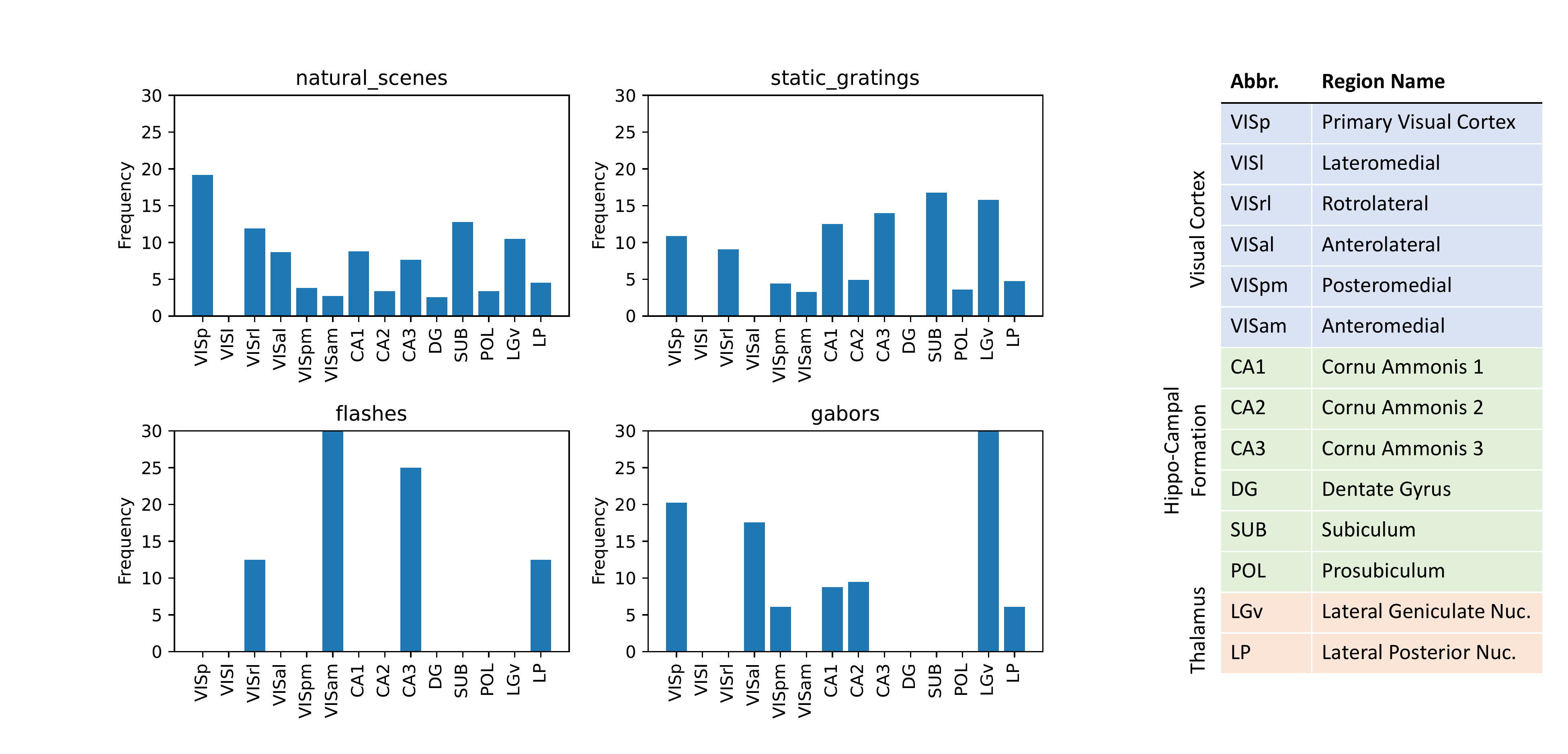}
    \caption{Proportional Frequency of appearance of each brain region in the hyper-edges of the estimated 4-order interaction structure.}
    \label{fig:barplot}
\end{figure}

\begin{figure}[t!]
    \centering
    \includegraphics[width = \textwidth]{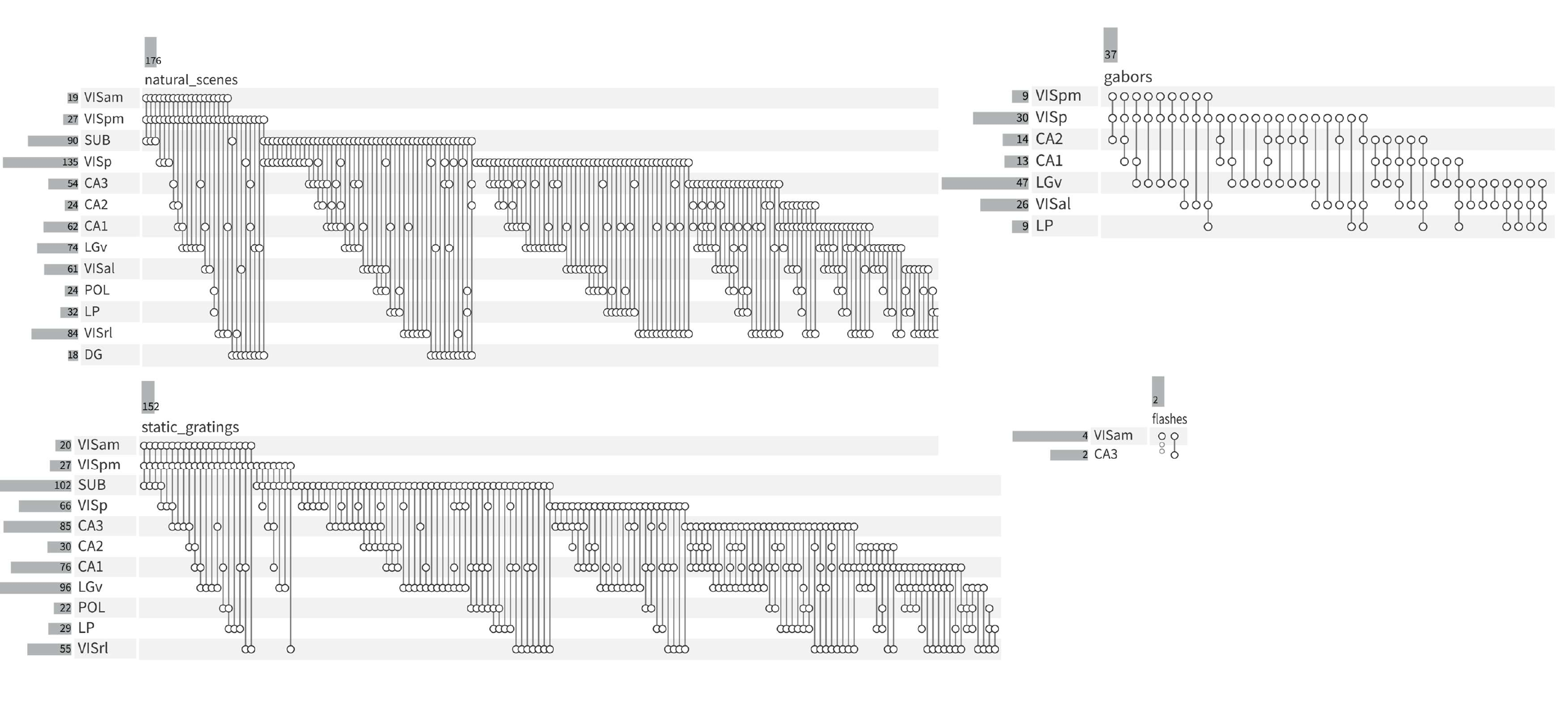}
    \caption{Estimated multi-body interactions between different regions of the mouse brain. For each vertical line, the small circles on it, represent the nodes that form the corresponding hyperedge.}
    \label{fig:paohvis}
\end{figure}

In Figure \ref{fig:barplot}, we show the proportional frequencies of occurrence of different brain regions in the $4$-body interaction hyper-edges. In the natural scenes scenario, VISp has the highest frequency of occurrence, followed by SUB, VISrl, VISal, LGv, CA1 and CA3. In the static grating scenario, SUB occur with the highest frequencies, followed by LGv, CA3 and  CA1. In the Gabor scenario, LGv occurs with the highest frequency, followed by VISp and VISal. Furthermore, many of the regions do not occur in any 4-body interactions, which can be explained by the fact that Gabor patches comprise of a gray screen with a small patch of gratings in a corner of the visual field of the mouse. Therefore Gabor patches is a less informative stimulus compared to natural scenes and static gratings, which can be attributed to the absence of several brain regions in the 4-body interactions. In the flashes scenario, VISam, CA3 appear with the highest frequencies, followed by VISrl and LP. Furthermore, we find that for natural scenes the regions in the visual cortex have relatively higher frequencies of occurence than the hippocampal formation and thalamus regions overall. Again, for the static gratings, the hippocampal formation have relatively higher frequency overall, compared to the other regions. In flashes, many of the brain regions do not occur in 4-body interactions, compared to the other three stimuli. This can be attributed to flashes comprising only of white or dark screen applied repeatedly to the mouse's visual field, thereby having less information present in the visual field. Although Gabor is also a less informative stimulus overall, in contrast to flashes, Gabor comprises a patch of informative gratings restricted to a small portion of the visual field. This can lead to a few brain regions having relatively greater frequency of higher order interactions, as seen in Figure \ref{fig:paohvis}.

%\textcolor{magenta}{any other pics?}
\section{Proofs of the Main Results}\label{mainprff}
The primary ingredient of the proof is the method of \textit{primal-dual witness} (see \cite{wainwright, wsharp}). 
Consider a primal solution $\widehat{\bm J}_r \in \mathbb{R}^{\bpk}$ along with an associated subgradient vector $\widehat{\bm z} \in \mathbb{R}^{\bpk}$ (which can be interpreted as a dual solution), such that the zero sub-gradient optimality conditions associated with the convex program (\ref{min}) are satisfied, which takes the  following form:
\begin{equation}\label{subgra1}
    \nabla \ell(\widehat{\bm J}_r)+\lambda \widehat{\bm z}=0,
\end{equation}
where the dual or subgradient vector $\widehat{\bm z} \in \mathbb{R}^{\bpk}$ must satisfy the properties
\begin{equation}\label{subgra2}
    \widehat{z}_{r,r_1,...,r_{k-1}}=\operatorname{sgn}(\widehat{J}_{r,r_1...r_{k-1}}) \quad \text{if } \widehat{J}_{r,r_1...r_{k-1}} \neq 0 \quad \text{and} \quad |\widehat{z}_{r,r_1,...,r_{k-1}}| \leq 1 \quad \text{otherwise.}
\end{equation}
By convexity, a pair $(\widehat{\bm J}_r, \widehat{\bm z}) \in \mathbb{R}^{\bpk} \times \mathbb{R}^{\bpk}$ is a primal-dual optimal solution to the convex program and its dual if and only if the two conditions (\ref{subgra1}) and (\ref{subgra2}) are satisfied. Note that the necessary and sufficient conditions that an optimal primal-dual pair correctly specifies the signed neighborhood of node $r$, are given by:
\begin{equation}\label{crt1}
    \mathrm{sgn}(\widehat{z}_{r, r_1,...,r_{k-1}})=\mathrm{sgn}(J_{r,r_1,...,r_{k-1}}) \quad \forall\{r_1,...,r_{k-1}\} \in \mathcal{N}(r) \quad \text { and}
\end{equation}
\begin{equation}\label{crt2}
    \widehat{J}_{r,r_1,...,r_{k-1}}=0 \quad \forall \{r_1,...,r_{k-1}\} \notin \mathcal{N}(r).
\end{equation}

Although the $\ell_{1}$-regularized logistic regression problem (\ref{min}) is convex, for $p \gg n$, it need not be strictly convex, and hence there may be multiple optimal solutions. In the following lemma, we provide sufficient conditions for shared sparsity among optimal solutions, as well as uniqueness of the optimal solution:

\begin{lem}\label{lem1}
Suppose that $\widehat{\bm J}$ is an optimal primal solution  with associated optimal dual vector $\widehat{\bm z}$ satisfying $\|\widehat{
\bm z}_{S^c}\|_\infty<1$. Then, any optimal primal solution $\widetilde{\bm J}$ must satisfy $\widetilde{\bm J}_{S^{c}}=0$. Moreover, if the Hessian sub-matrix $\left[\nabla^{2} \ell(\widehat{\bm J})\right]_{SS}$ is positive definite, then $\widehat{\bm J}$ is the unique optimal solution.
\end{lem}

Lemma \ref{lem1} is proved in Appendix \ref{lem1apx}. It helps us to construct a primal-dual witness $(\hat{\bm J},\hat{\bm z})$ in the following steps:

\begin{enumerate}
    \item Define:
    \begin{equation}\label{gu1}
        \widehat{\bm J}_{S}=\underset{(\bm J_{S}, \boldsymbol{0}) \in \mathbb{R}^{\binom{p-1}{k-1}}}{\arg \min }\{\ell(\bm J ; \mathfrak{X}^{n})+\lambda\|\bm J_{S}\|_{1}\}
    \end{equation}
    and set $\widehat{\bm z}_{S}=\operatorname{sgn}(\widehat{\bm J}_S)$.
\vspace{0.1in}

    \item Next, set $\widehat{\bm J}_{S^c} = \boldsymbol{0}$
    \vspace{0.1in}

    \item We then obtain $\widehat{\bm z}_{S^c}$ by substituting in \eqref{subgra1} the values of $\widehat{\bm J}$ and $\widehat{\bm z}_S$.
    \vspace{0.1in}

    \item Finally, we show that under the assumptions of Theorem \ref{thm1}, conditions \eqref{crt1} and \eqref{subgra2} are satisfied with high probability.
\end{enumerate}

It will become evident from the proof in step (4), that $\left\|\widehat{\bm z}_{S^{c}}\right\|_{\infty}<1$ with high probability. Also, we will prove that the Hessian sub-matrix $\left[\nabla^2 \ell(\hat{\bm J})\right]_{SS}$ is strictly positive definite with high probability which will enable us to conclude in view of Lemma \ref{lem1}, that the primal solution $\widehat{\bm J}$ is unique.

\subsection{Some Technical Lemmas}
The first step of our analysis is to show the consistency of $\widehat{\bm J}$ under certain conditions on $(n,p,d,k,\lambda)$ and assumptions on the sample Fisher and covariance matrices: 
$$\widehat{\bm Q} := -\widehat{\mathbb{E}}[\nabla^2\log \mathbb{P}_{\bm J}(X_r|\bm X_{\backslash r})] := \frac{1}{n} \sum_{i=1}^n \eta_r(\bm x^{(i)}; \bm J) \bm x_{\cdot r}^{(i)} \left[\bm x_{\cdot r}^{(i)}\right]^\top\quad \text{and}\quad \hat{\boldsymbol{\Sigma}} := \frac{1}{n}\sum_{i=1}^n \bm x_{\cdot r}^{(i)} (\bm x_{\cdot r}^{(i)})^\top~.$$ Specifically, we assume for sometime, that Assumptions \eqref{a1} and \eqref{a2} hold on the sample Fisher matrix $\widehat{\bm Q}$ and the sample covariance matrix $\widehat{\boldsymbol{\Sigma}}$. We begin with some technical results.
We start by observing that the zero subgradient condition can be written as:
\begin{equation}\label{thrt}
    \nabla \ell(\widehat{\bm J} ; \mathfrak{X}^{n})-\nabla \ell(\bm J ; \mathfrak{X}^{n})=\bm W-\lambda\widehat{\bm z},
\end{equation}
with $\bm W := -\nabla \ell(\bm J ; \mathfrak{X}^{n})$ for the $\bpk$-dimensional score function,
$$
\bm W:= \frac{k!}{n} \sum_{i=1}^{n}\bm x_{\cdot r}^{(i)}\left\{x_{r}^{(i)}-\frac{e^{km_r(\bm x^{(i)})}-e^{-km_r(\bm x^{(i)})}}{e^{km_r(\bm x^{(i)})}+e^{-km_r(\bm x^{(i)})}}\right\} .
$$

It is easy to see that $\mathbb{E}_{\bm J}\left(\bm W\right)=0$. Next, applying the mean-value theorem to \eqref{thrt}, we have:
\begin{equation}\label{subgra3}
    \nabla^{2} \ell (\bm J; \mathfrak{X}^{n})[\widehat{\bm J}-\bm J]=\bm W-\lambda \widehat{\bm z}-\bm R,
\end{equation}
where the remainder term is given by:
\begin{equation}\label{remterm8}
    \quad R_{j} :=[\nabla^{2} \ell(\bar{\bm J}^{(j)} ; \mathfrak{X}^{n})-\nabla^{2} \ell(\bm J ; \mathfrak{X}^{n})]_{j*}(\widehat{\bm J}- \bm J)
\end{equation}
with $\bar{\bm J}^{(j)}$ being a parameter vector on the line segment joining $\bm J$ and $\hat{\bm J}$, and $[\cdot]_{j*}$ denoting the $j^{\mathrm{th}}$ row of the matrix. Below, we provide a probabilistic upper bound on the term $\bm W$.

\begin{lem}\label{lem2}
For the mutual incoherence parameter $\alpha \in(0,1]$ in Assumption \ref{a27}, we have
$$
\mathbb{P}\left(\frac{2-\alpha}{\lambda}\left\|\bm W\right\|_{\infty} \geq \frac{\alpha}{4}\right) \leq 2 \exp \left(-\frac{n\alpha^2 \lambda^2}{128(2-\alpha)^2 (k!)^2} +\log {p-1 \choose k-1}\right).
$$
\end{lem}

 Lemma \ref{lem2} is proved in Appendix \ref{aplem2}. The following result establishes $\ell^2$-consistency of $\widehat{\bm J}_S$ towards $\bm J_S$. 

\begin{lem}\label{l2con}
If $\lambda \leq \frac{C_{\min }^2}{40 d D_{\max} (k!)^3}$ and $\|\bm W\|_{\infty} \leq \lambda / 4$, then
$$
\|\widehat{\bm J}_S-\bm J_S\|_2 \leq \frac{5}{2 C_{\min }} \lambda 
 \sqrt{d}.
$$
\end{lem}
Lemma \ref{l2con} is proved in Appendix \ref{l2conpr}. Finally, we control the remainder term \eqref{remterm8}.
\begin{lem}\label{l337}
If $\lambda d \leq \frac{C_{\min }^2}{200 (k!)^3 D_{\max }} \cdot \frac{\alpha}{2-\alpha}$ and $\left\|\bm W\right\|_{\infty} \leq \lambda / 4$, then
$$
\frac{\left\|\bm R\right\|_{\infty}}{\lambda} \leq50 (k!)^3 \lambda d\frac{D_{\max}}{C_{\min}^2} \leq \frac{\alpha}{4(2-\alpha)}$$
\end{lem}
Lemma \ref{l337} is proved in Appendix \ref{apl337}.
 
\subsection{Recovery Under Sample Assumptions}\label{sec:samplerec}
To begin with, we show that Theorem \ref{thm1} holds if we assume conditions \eqref{a1} and \eqref{a2} on the sample Fisher matrix $\widehat{\bm Q}$ and the sample covariance matrix $\widehat{\boldsymbol{\Sigma}}$. Moreover, this is true under the slightly weaker assumption $n> L (k!)^8 d^2 \log \bpk$ on the sample size.

To begin with, choose $\lambda = 16 k! [(2-\alpha)/\alpha]\sqrt{n^{-1} \log \bpk}$. It then follows from Lemma \ref{lem2} that with probability $1-\bpk^{-1}$, we have:
$$
\|\bm W\|_{\infty} \leq \frac{\alpha\lambda }{4(2-\alpha)} \leq \frac{\lambda}{4}~.
$$
Next, we verify the remaining two conditions in the hypotheses of Lemma \ref{l2con} and Lemma \ref{l337}. For this, note that the upper bound on $\lambda$ specified in the hypothesis of Lemma \ref{l2con} is implied by the the upper bound on $\lambda$ specified in the hypothesis of Lemma \ref{l337}, and hence, it is sufficient to verify the latter only. Towards this, note that:
    $$\lambda d \leq \frac{C_{\min }^2}{200 (k!)^3 D_{\max }} \cdot \frac{\alpha}{2-\alpha}\iff n \ge 3200^2 (k!)^8d^2 \left(\frac{2-\alpha}{\alpha}\right)^4 \frac{D_{\max}^2}{C_{\min}^4} \log \bpk$$
and hence, we may choose $L > 3200^2[(2-\alpha)/\alpha]^4 D_{\max}^2/C_{\min}^4$, so that the hypotheses of both Lemma \ref{l2con} and Lemma \ref{l337} are satisfied.

Next, let $\widehat{\bm J}_{S}$ be the minimizer of the partial penalized likelihood \eqref{gu1}. Set $\widehat{\bm z}_{S}=\text{sign}(\widehat{\bm J}_{S})$ and $\widehat{\bm J}_{S^{c}}=0$. Then we can find a $\widehat{\bm z}_{S^c}$ that makes $\widehat{\bm J}$ a solution to the original convex program, thereby satisfying \eqref{subgra3}, which can be re-written in block form as:
\begin{equation}\label{conblock}
    \begin{aligned}
        \widehat{\bm Q}_{S^c S}[\widehat{\bm J}_S- \bm J_S]&= \bm W_{S^c}-\lambda \widehat{\bm z}_{S^c} - \bm R_{S^c},\\
        \widehat{\bm Q}_{S S}[\widehat{\bm J}_S-\bm J_S]&= \bm W_S -\lambda \widehat{\bm z}_S - \bm R_S,
    \end{aligned}
\end{equation}
It thus follows from \eqref{conblock} that:

\begin{equation*}
    \quad \widehat{\bm Q}_{S^c S} \widehat{\bm Q}_{S S}^{-1}\left[\bm W_S-\lambda \widehat{\bm z}_S- \bm R_S\right]=\bm W_{S^c} -\lambda \widehat{\bm z}_{S^c} - \bm R_{S^c}
\end{equation*}
which, on rearrangement, yields the following:
\begin{equation}\label{rearr17}
    \bm W_{S^c}- \bm R_{S^c} - \widehat{\bm Q}_{S^c S} \widehat{\bm Q}_{S S}^{-1}\left[\bm W_S- \bm R_S\right]+\lambda \widehat{\bm Q}_{S^c S}\widehat{\bm Q}_{S S}^{-1} \widehat{\bm z}_S=\lambda \widehat{\bm z}_{S^c}
\end{equation}

We now aim to show that $\widehat{\bm J}$ is a unique solution to the convex program. Towards this, note that by \eqref{rearr17}, we have:
\begin{eqnarray*}
    \|\widehat{\bm z}_{S^c}\|_{\infty} &\leqslant& \|\widehat{\bm Q}_{S^c S} \widehat{\bm Q}_{S S}^{-1}\|_{\infty}\left[\frac{\|\bm W_S\|_{\infty}}{\lambda}+\frac{\|\bm R_S\|_{\infty}}{\lambda}+1\right]
    +\frac{\|\bm R_{S^c}\|_{\infty}}{\lambda}+\frac{\|\bm W_{S^c}\|_{\infty}}{\lambda} \\
    &\leqslant& 1-\alpha + (2-\alpha)\left[\frac{\|\bm R\|_{\infty}}{\lambda}+\frac{\|\bm W\|_{\infty}}{\lambda}\right]\\&\leqslant& 1-\alpha +\frac{\alpha}{4}+\frac{\alpha}{4}=1-\frac{\alpha}{2}<1
\end{eqnarray*}
with probability converging to one, by Lemma \ref{l337}. By Lemma \ref{lem1} and Assumption \eqref{a1} on the sample Fisher matrix, we can thus conclude that $\widehat{\bm J}$ is the unique solution to the $\ell_1$-regularized logistic regression. 

Next, we show that $\widehat{\bm J}_S$ defined by (\ref{gu1}) satisfies $\mathrm{sgn}(\widehat{\bm J}_S)=\mathrm{sgn}(\bm J_S)$. In order to do so, it suffices to show that
$$
\|\widehat{\bm J}_S- \bm J_S\|_{\infty} \leq \frac{J_{\min }}{2}
$$
where $J_{\min }:=\min _{\bm e \in E}|J_{\bm e}|$. From Lemma \ref{l2con}, we have $\| \widehat{\bm J}_S- \bm J_S \|_2 \leq \frac{5}{2C_{\min }} \sqrt{d} \lambda$ so that
$$
\begin{aligned}
    \frac{2}{J_{\min }}\|\widehat{\bm J}_S- \bm J_S\|_{\infty} & \leq  \frac{2}{J_{\min }}\|\widehat{\bm J}_S- \bm J_S\|_2 \\
    & \leq \frac{2}{J_{\min }}\cdot  \frac{5}{2C_{\min }} \sqrt{d} \lambda,
\end{aligned}
$$
which is less than $1$ as long as $J_{\min } \geq \frac{5}{C_{\min }} \sqrt{d} \lambda$. This concludes the proof of Theorem \ref{thm1} under the assumptions on the sample Fisher and sample covariance matrices.

\subsection{Transferring from Sample to Population Assumptions}
In this section, we show that Assumptions \ref{a1} and \ref{a2} on the population Fisher and covariance matrices actually imply analogous bounds on their sample versions. This will then enable us to conclude Theorem \ref{thm1} in view of Section \ref{sec:samplerec}. We start with a result guaranteeing high probability occurrence of the eigenvalue bounds \eqref{a1} for the sample Fisher and covariance matrices, under the corresponding population assumptions \eqref{a1}.

% In the proof of the following lemmas, a bound on the terms of Hessian matrix will be extensively used. By definition of the matrices $Q^n$ and $Q^*$ (as in (\ref{qstar})), the $(j, k)$ th element of the difference matrix $Q^n-Q^*$ can be written as an i.i.d. sum of the form $Z_{s t}=\frac{1}{n} \sum_{i=1}^n Z_{s t}^{(i)}$ where each $Z_{s t}^{(i)}$ is zero-mean and bounded (in particular, $\left|Z_{s t}^{(i)}\right| \leq 2 k^2$ ). Similar to the proof of Lemma \ref{lem2}, by the Azuma-Hoeffding inequality, for any indices $j, k=1, \ldots, d$ and for any $\varepsilon>0$, we have
% \begin{equation}\label{hoeffding}
%     \mathbb{P}\left[\left(Z_{s t}\right)^2 \geq \varepsilon^2\right]=\mathbb{P}\left[\left|\frac{1}{n} \sum_{i=1}^n Z_{s t}^{(i)}\right| \geq \varepsilon\right] \leq 2 \exp \left(-\frac{\varepsilon^2 n}{2k^4}\right).
% \end{equation}

\begin{lem}\label{lem5}
Under Assumption \eqref{a1}, for every $\delta>0$ and constants $A, B$, we have:
\begin{equation}\label{lem5eqn1}
    \mathbb{P}\left[\Lambda_{\max}\left[\frac{1}{n}\sum_{i=1}^n \bm x_{\cdot r}^{(i)}(\bm x_{\cdot r}^{(i)})^\top\right] \geqslant D_{\max }+\delta\right] \leqslant 2 \exp \left(-A \frac{n\delta^2}{d^2 (k!)^4}+B \log d\right),
\end{equation}
\begin{equation}\label{lem5eqn2}
    \mathbb{P}\left[\Lambda_{\min }\left(\widehat{\bm Q}_{S S}\right) \leqslant C_{\min }-\delta\right] \leqslant 2 \exp \left(-A \frac{n \delta^2}{d^2 (k!)^4}+B \log d\right),
\end{equation}
\end{lem}

Lemma \ref{lem5} is proved in Appendix \ref{aplem5}. Now we state the analogous result for the incoherence condition \eqref{a2}, which guarantees sample incoherence from population incoherence.

\begin{lem}\label{sampinc}
Suppose that the population Fisher matrix satisfies the incoherence condition (\ref{a2}) with parameter $\alpha \in(0,1]$. Then the sample Fisher matrix satisfies:
\begin{equation}\label{lem6eqn}
    \quad \mathbb{P}\left[\vertiii{\widehat{\bm Q}_{S^c S} \widehat{\bm Q}_{S S}^{-1}}_{\infty} \geq 1-\frac{\alpha}{2}\right] \leq M \exp \left(-K_1 \frac{n}{d^3 (k!)^4}+ K_2 \log d + \log \left[{p-1 \choose k-1}-d\right]\right)
\end{equation}
for some positive constants $M,K_1$ and $K_2$.
\end{lem}

Lemma \ref{sampinc} is proved in Appendix \ref{apsampinc}. With Lemmas \ref{lem5} and \ref{sampinc} in hand, we are now ready to prove Theorem \ref{thm1}. Towards this, define
$G_{n,\delta,\alpha} := G_{1,n,\delta}\bigcup G_{2,n,\delta}\bigcup G_{3,n,\alpha}$, where:

$$G_{1,n,\delta} := \left\{\Lambda_{\max}\left[\frac{1}{n}\sum_{i=1}^n \bm x_{\cdot r}^{(i)}(\bm x_{\cdot r}^{(i)})^\top\right] \geqslant D_{\max }+\delta\right\}~,$$

$$G_{2,n,\delta} := \left\{\Lambda_{\min }\left(\widehat{\bm Q}_{S S}\right) \leqslant C_{\min }-\delta\right\}~\quad\text{and}$$

$$G_{3,n,\alpha} :=\left\{\vertiii{\widehat{\bm Q}_{S^c S} \widehat{\bm Q}_{S S}^{-1}}_{\infty} \geq 1-\frac{\alpha}{2}\right\}~.$$
In view of Lemmas \ref{lem5} and \ref{sampinc}, we have:
$$\p\left(G_{n,\delta,\alpha}\right) \le M \exp \left(-K_1(\delta) \frac{n}{d^3 (k!)^4}+ K_2 \log d + \log \left[{p-1 \choose k-1}-d\right]\right)$$
for some constants $M,K_1,K_2>0$. On the other hand, on the event $G_{n,\delta,\alpha}^c$, the conclusions of Section \ref{sec:samplerec} hold. This completes the proof of Theorem \ref{thm1}.

\section{Discussion}\label{discuss_sec}
In this paper, we showed that a simple node-wise $\ell_1$-regularized logistic regression technique inspired by \cite{wainwright} can be used to consistently recover the tensor structure of a $k$-spin Ising model. We worked in the high-dimensional setting allowing both the dimension $p$ and the interaction factor $k$ of the model to grow with the number of samples $n$. Our result shows that consistent tensor recovery is possible for sample sizes $n = \Omega((k!)^8 d^3 \log \bpk)$, where $d$ denotes the maximum degree of the Ising hypergraph. Our theoretical results are supported by the two simulation settings we considered in Section \ref{simstudy}, where in each case, the hyperedge recovery rate is seen to approach $1$ with increase in sample size. We also applied our method on a real-life neurobiological dataset involving electro-psychological recordings from the mouse brain, and revealed higher-order neural interactions between the different regions of the mouse brain.   

This work also leaves some potentially interesting directions for future research. One of these areas is to prove the minimax optimality of the $\min\{1,\sqrt{|E|/p}\}$ rate of learning tensor Ising models, as discussed in Remark \ref{lucrem}. Another potential goal may be to improve the computational complexity of the algorithm considered in this paper, by possibly adapting the methods considered in \cite{bresler, bresler1, bresler2} in the tensor setting. A third direction for future work may be to consider the dependent sampling scenario (for example, data coming from a time series). An interesting question under this setting, is that whether consistent structure learning is possible under some weak dependence frameworks, such as $\rho$-mixing or strong mixing.

\section{Acknowledgement}
S.M. was supported by  by the National University of Singapore Start-Up Grant R-155-000-233-133, 2021. The authors are grateful to Luc Devroye for several helpful discussions on the Ising model learning problem mentioned in Remark \ref{lucrem}.

\appendix
\section{Proof of Lemma \ref{lem1}}\label{lem1apx}
Note that the penalized problem (\ref{min}) can be equivalently written as a constrained optimization problem over the ball $\{\|\bm J\|_{1} \leq C\left(\lambda\right)\}$, for some constant $C\left(\lambda\right)<+\infty$ (by Lagrange duality). The primal can be formulated as:
\begin{equation}\label{primal}\tag{$\mathbb{P}$}
    \begin{aligned}
    &\text{min} \quad \ell(\bm J),\\
    &\text{s.t.} \quad \|\bm J\|_1 - C\leqslant 0,\\
    & \bm J\in \mathbf{R}^{(p-1)^{k-1}}.
    \end{aligned}
\end{equation}
whose dual can be written as:
\begin{equation}\label{dual}\tag{$\mathbb{D}$}
    \begin{aligned}
        &\text{max}\quad \inf\limits_{\bm J\in\mathbf{R}^{(p-1)^{k-1}}} \ell(\bm J)+\lambda(\|\bm J\|_1-C),\\
        &\text{s.t.} \quad \lambda \geqslant 0.\\
    \end{aligned}
\end{equation}
If there exists a feasible $\widehat{\bm J}$ that optimises \ref{primal}, by strong duality we have 
\begin{equation}\label{slack}
    \begin{aligned}
        &\lambda(\|\widehat{\bm J}\|_1-C)=0,\\
        \text{and}\quad &\nabla l(\widehat{\bm J}) + \lambda \widehat{\bm z}=0,
    \end{aligned}
\end{equation}
Note that by definition of the $\ell_{1}$-subdifferential, the subgradient vector $\widehat{\bm z}$ can be expressed as a convex combination of sign vectors of the form 
\begin{equation}
    \widehat{\bm z}=\sum\limits_{\bm v \in\{-1,1\}^{(p-1)^{k-1}}} \alpha_{\bm v} \bm v, 
\end{equation}
where the weights $\alpha_{\bm v}$ form a probability vector. We consider an alternative formulation of a pair of primal-dual problem, given by: 
\begin{equation}\label{primal1}\tag{$\mathbb{P}'$}
    \begin{aligned}
    &\text{min} \quad \ell(\bm J),\\
    &\text{s.t.} \quad \langle \bm v, \bm J \rangle - C\leqslant 0, \quad \forall \bm v\in\{-1,+1\}^{(p-1)^{k-1}},\\
    & \bm J\in \mathbf{R}^{(p-1)^{k-1}}
    \end{aligned}
\end{equation}
and
\begin{equation}\label{dual1}\tag{$\mathbb{D}'$}
    \begin{aligned}
        &\text{max}\quad \inf\limits_{\bm J\in\mathbf{R}^{(p-1)^{k-1}}} \ell(\bm J)+\lambda \sum\limits_{\bm v \in\{-1,+1\}^{(p-1)^{k-1}}} \alpha_{\bm v}(\langle \bm v, \bm J\rangle-C),\\
        &\text{s.t.} \quad \lambda \geqslant 0,~\alpha_{\bm v} \geqslant 0.\\
    \end{aligned}
\end{equation}
We notice that \ref{primal} is equivalent to \ref{primal1}. Then, any other optimal primal solution $\widetilde{\bm J}$ also satisfies the complementary slackness conditions and the zero subgradient optimality condition (as in \ref{slack}):

\begin{equation}\label{slack1}
    \begin{aligned}
        &\lambda \sum\limits_{\bm v \in\{-1,1\}^{(p-1)^{k-1}}} \alpha_{\bm v}(\langle \bm v,\widetilde{\bm J}\rangle-C)=0,\\
        \text{and}\quad &\nabla \ell (\widetilde{\bm J}) + \lambda \widehat{\bm z}=0.
    \end{aligned}
\end{equation}

As $\lambda > 0$, the slackness conditions in (\ref{slack}) and (\ref{slack1}) imply that $\langle\widehat{\bm z}, \widetilde{\bm J}\rangle=C=\|\widetilde{\bm J}\|_{1}$ which is impossible if $\widetilde{J}_{\bm e} \neq 0$ for some index $\bm e$ for which $|\widehat{z}_{\bm e}|<1$. Since $\|\widehat{\bm z}_{S^c}\|_\infty < 1$, it follows that $\widetilde{\bm J}_{S^{c}}=0$ for all optimal primal solutions $\widetilde{\bm J}$. This proves the first part of Lemma \ref{lem1}.

For proving the second part, observe that since all optimal solutions satisfy $\widetilde{\bm J}_{S^{c}}=0$, we may restrict our optimization problem to this set of constraints. If the principal submatrix $\left[\nabla^{2} \ell(\widehat{\bm J})\right]_{SS}$ of the Hessian is positive definite, then this sub-problem is strictly convex, thereby guaranteeing a unique optimal solution. This completes the proof of Lemma \ref{lem1}.

\section{Proof of Lemma \ref{lem2}}\label{aplem2}
The coordinates of $\bm W$ can be written as $W_{\bm s} := \frac{k!}{n} \sum_{i=1}^n Z_{\bm s}^{(i)}$, where
$$
Z_{\bm s}^{(i)}:=x_{\bm s}^{(i)}\left\{x_r^{(i)}-\mathbb{P}_{\bm J}\left[x_r^{(i)}=1 \mid x_{\backslash r}^{(i)}\right]+\mathbb{P}_{\bm J}\left[x_r^{(i)}=-1 \mid x_{\backslash r}^{(i)}\right]\right\}
$$
with $x_{\bm s} := x_{s_1}\ldots x_{s_{k-1}}$. Note that under $\p_{\bm J}$, the random variables $\{Z_{\bm s}^{(i)}\}_{1\le i\le n}$ are i.i.d. with mean zero, and bounded (by $2$).
By Hoeffding inequality, we thus have:

$$\mathbb{P}\left(\left|W_{\bm s}\right|>\delta\right)\leq 2 \exp \left(-\frac{n \delta^2}{8(k!)^2}\right)$$ Setting $\delta := \alpha \lambda /4(2-\alpha)$, we get:

$$\mathbb{P}\left(\frac{2-\alpha}{\lambda}\left|W_{\bm s}\right|>\frac{\alpha}{4}\right) \leq 2 \exp \left(-\frac{n\alpha^2 \lambda^2}{128(2-\alpha)^2 (k!)^2}\right) .
$$
A union bound over the indices $\bm s \in T_r$ now gives:
$$
\mathbb{P}\left(\frac{2-\alpha}{\lambda}\left\|W_{\bm s}\right\|_{\infty}>\frac{\alpha}{4}\right) \leq 2 \exp \left(-\frac{n\alpha^2 \lambda^2}{128(2-\alpha)^2 (k!)^2} +\log {p-1 \choose k-1}\right)
$$
which completes the proof of Lemma \ref{lem2}.

\section{Proof of Lemma \ref{l2con}}\label{l2conpr}

Define $G: \mathbb{R}^d \rightarrow \mathbb{R}$ by:
\begin{equation}\label{gu}
    G\left(\bm t\right):=\ell\left(\bm J_S +\bm t ; \mathfrak{X}^n\right)-\ell\left(\bm J_S ; \mathfrak{X}^n\right)+\lambda\left(\left\|\bm J_S+\bm t\right\|-\left\|\bm J_S\right\|\right).
\end{equation}
Note that $G$ is a convex function. It follows from \eqref{gu1} that $\widehat{\bm t}=\widehat{\bm J}_S-\bm J_S$ minimizes $G$. Also, since $G(\boldsymbol{0})=0$, we have $G(\widehat{\bm t}) \leq 0$. We claim that if we can show that $G(\bm t) > 0$ for all $\bm t\in \mathbb{R}^d$ with $\|\bm t\|_2 = B$ for some $B >0$, then $\|\widehat{\bm t}\|_2 \le B$. To see this claim, note that if $\widehat{\bm t}$ lied outside the ball of radius $B$, then the vector $\|\delta \widehat{\bm t}\|_2 = B$ for some appropriately chosen $\delta \in(0,1)$. By convexity,
$$
G(\delta \widehat{\bm t}) \leq \delta G(\widehat{\bm t})+(1-\delta) G(\boldsymbol{0}) \leq 0,
$$
contradicting the assumed strict positivity of $G$ on the boundary. 

In view of all these, it is thus enough to prove that $G(\bm t) > 0$ for all $\bm t$ satisfying $\|\bm t\|_2 = B := M \lambda \sqrt{d}$ where $M>0$ is to be chosen later in the proof. Towards this, let $\bm t \in \mathbb{R}^d$ satisfy $\|\bm t\|_2=B$. Then, by a Taylor series expansion, we have:
\begin{equation}\label{gexp}
    G(\bm t)= \bm W_S^\top \bm t + \bm t^\top\left[\nabla^2 \ell\left(\bm J_S+\alpha \bm t\right)\right] \bm t+\lambda\left(\left\|\bm J_S+\bm t\right\|-\left\|\bm J_S\right\|\right)
\end{equation}
for some $\alpha \in[0,1]$, where $\bm W := \nabla \ell(\bm J; \mathfrak{X}^n)$. For the first term, we have the bound
\begin{equation}\label{gexp7}
    \quad|\bm W_S^\top \bm t| \leq\|\bm W_S\|_{\infty}\|\bm t\|_1 \leq\|\bm W_S\|_{\infty} \sqrt{d}\|\bm t\|_2 \leq \frac{M\lambda^2 d}{4},
\end{equation}
since by assumption, $\|\bm W_S\|_{\infty} \leq \frac{\lambda}{4}$ by assumption.
Also, an application of the triangle inequality gives:
$$
\lambda \left(\|\bm J_S+\bm t\|_1-\|\bm J_S\|_1\right) \geq-\lambda\|\bm t\|_1 \ge -\lambda \sqrt{d} \|\bm t\|_2 = - M\lambda^2 d.
$$

Finally, we analyze the middle quadratic form. Towards this, we note that:

\begin{eqnarray*}
   \zeta &:=& \Lambda_{\min}\left(\nabla^2 \ell\left(\bm J_S +\alpha \bm t ; \mathfrak{X}^n\right)\right) \\
&\geq& \min _{\alpha \in[0,1]} \Lambda_{\min }\left(\nabla^2 \ell\left(\bm J_S+\alpha \bm t ; \mathfrak{X}^n\right)\right) \\
&=& \min _{\alpha \in[0,1]} \Lambda_{\min }\left[\frac{1}{n} \sum_{i=1}^n \eta\left(\bm x^{(i)} ; \bm J_S +\alpha \bm t\right) \bm x_S^{(i)}\left(\bm x_S^{(i)}\right)^\top\right] 
\end{eqnarray*}
where $\eta := \eta_r$. We therefore have by a Taylor series expansion,

\begin{eqnarray*}
\zeta &\geq & \Lambda_{\min}\left[\frac{1}{n} \sum_{i=1}^n \eta(\bm x^{(i)} ; \bm J_S) \bm x_S^{(i)}(\bm x_S^{(i)})^\top\right] \\
&-& \max_{\alpha \in[0,1]}\vertiii{\frac{1}{n} \sum_{i=1}^n \eta^{\prime}(\bm x^{(i)} ; \bm J_S+\alpha \bm t)(\bm t^\top \bm x_S^{(i)}) \bm x_S^{(i)}(\bm x_S^{(i)})^\top}_2\\
&=& \Lambda_{\min }(\widehat{\bm Q}_{SS})-\max _{\alpha \in[0,1]}\vertiii{\frac{1}{n} \sum_{i=1}^n \eta^{\prime}(\bm x^{(i)} ; \bm J_S+\alpha \bm t)(\langle \bm t, \bm x_S^{(i)}\rangle) \bm x_S^{(i)}(\bm x_S^{(i)})^\top}_2\\
&\geq & C_{\min }-\max _{\alpha \in[0,1]}\vertiii{\underbrace{\frac{1}{n} \sum_{i=1}^n \eta^{\prime}(\bm x^{(i)} ; \bm J_S+\alpha \bm t)(\langle \bm t, \bm x_S^{(i)}\rangle) \bm x_S^{(i)}(\bm x_S^{(i)})^\top}_{\bm M(\alpha)}}_2 .
\end{eqnarray*}
Moving forward, the aim is thus to control the spectral norms of the matrices $\bm M(\alpha)$, for $\alpha \in[0,1]$. Towards this, note that for any fixed $\alpha \in[0,1]$ and $\bm y$ lying on the boundary of the unit ball in $\mathbb{R}^d$, we have:

\begin{eqnarray*}
\langle \bm y, \bm M (\alpha) \bm y\rangle &=& \frac{1}{n} \sum_{i=1}^n \eta^{\prime}(\bm x^{(i)} ; \bm J_S+\alpha \bm t) \langle \bm t, \bm x_S^{(i)}\rangle \langle \bm x_S^{(i)}, \bm y\rangle^2 \\
& \leq& \frac{1}{n} \sum_{i=1}^n |\eta^{\prime}(\bm x^{(i)} ; \bm J_S+\alpha \bm t)| |\langle \bm t, \bm x_S^{(i)}\rangle| \langle \bm x_S^{(i)}, \bm y\rangle^2 .
\end{eqnarray*}
Now note that $|\eta^{\prime}(\bm x^{(i)} ; \bm J_S+\alpha \bm t)| \leq 8(k!)^3$, and
$$
|\langle \bm t, \bm x_S^{(i)}\rangle| \leq \sqrt{d}\|\bm t\|_2=M \lambda d .
$$
Also, note that by our sample assumptions,
$$
\frac{1}{n} \sum_{i=1}^n \langle \bm x_S^{(i)}, \bm y\rangle ^2 \leq \vertiii{\frac{1}{n} \sum_{i=1}^n \bm x_S^{(i)}(\bm x_S^{(i)})^T}_2 \leq D_{\max }
$$
It follows from the above inequalities, that
$$
\max _{\alpha \in[0,1]}\vertiii{\bm M (\alpha)}_2 \leq 8 D_{\max } M \lambda d  (k!)^3 \leq C_{\min } / 2
$$
as long as $\lambda \le C_{\min}/[16D_{\max}M d(k!)^3]$, which we will verify soon, after we specify the constant $M$.
Under this condition, we have shown that
\begin{equation}\label{gexp8}
    \zeta :=\Lambda_{\min }\left(\nabla^2 \ell\left(\bm J_S +\alpha \bm t ; \mathfrak{X}^n\right)\right) \geq C_{\min } / 2
\end{equation}

Finally, it follows from the \eqref{gexp}, \eqref{gexp7} and \eqref{gexp8} that:

$$
G\left(\bm t\right) \geq M \lambda^2 d \left(-\frac{1}{4} -1+\frac{M C_{\min}}{2}\right)
$$
the right-hand side being strictly positive if $M > 5/(2C_{\min})$.  Therefore, under the assumption:
$$
\lambda \leq \frac{C_{\min }}{16 D_{\max } M d (k!)^3} < \frac{C_{\min }^2}{40 D_{\max } d (k!)^3},
$$
we can conclude that:
$$
\|\widehat{\bm t}\|_2 \leq M \lambda \sqrt{d}=\frac{5}{2 C_{\min }} \lambda \sqrt{d}$$
thereby completing the proof of Lemma \ref{l2con}.

\section{Proof of Lemma \ref{l337}}\label{apl337}
To begin with, note that for every $j \in\left[{p-1 \choose k-1}\right]$, we have:
$$
\begin{aligned}
R_j &=[\nabla^2 \ell(\bar{\bm J}^{(j)} ; \mathfrak{X}^n)-\nabla^2 \ell(\bm J ; \mathfrak{X}^n)]_{j*} [\widehat{\bm J}- \bm J] \\
&=-\frac{1}{n} \sum_{i=1}^n[\eta(\bm x^{(i)} ; \bar{\bm J}^{(j)})-\eta(\bm x^{(i)} ; \bm J)][\bm x_{\cdot r}^{(i)} \bm (\bm x_{\cdot r}^{(i)})^T]_{j*}[\widehat{\bm J}-\bm J],
\end{aligned}
$$
for some point $\bar{\bm J}^{(j)}=t_j \widehat{\bm J}+\left(1-t_j\right) \bm J$. Setting $$h(t)=\frac{4(k!)^2 \exp (2k! t)}{( \exp (2k! t)+1)^2}~,$$ observe that $\eta(\bm x ; J)=h(x_r m_r(\bm x)/(k-1)!) = h(x_r \bm J^\top \bm x_{\cdot r})$. We thus have:

$$
\begin{aligned}
-R_j &= \frac{1}{n} \left[h(x_r^{(i)} (\widetilde{\bm J}^{(j)})^\top \bm x_{\cdot r}^{(i)}) - h(x_r^{(i)}\bm J^\top \bm x_{\cdot r}^{(i)})\right][\bm x_{\cdot r}^{(i)}  (\bm x_{\cdot r}^{(i)})^\top]_{j*}[\widehat{\bm J}-\bm J]\\&=\frac{1}{n} \sum_{i=1}^n h^{\prime}(x_r^{(i)} (\widetilde{\bm J}^{(j)})^\top \bm x_{\cdot r}^{(i)}) x_r^{(i)}(\bm x_{\cdot r}^{(i)})^T[\bar{\bm J}^{(j)}- \bm J][x_{\cdot r, j}^{(i)}(\bm x_{\cdot r}^{(i)})^\top][\widehat{\bm J}-\bm J]\} \\
&=\frac{1}{n} \sum_{i=1}^n\{h^{\prime}(x_r^{(i)}(\widetilde{\bm J}^{(j)})^\top \bm x_{\cdot r}^{(i)})  x_r^{(i)}x_{\cdot r, j}^{(i)}\}\{[\bar{\bm J}^{(j)}-\bm J]^\top \bm x_{\cdot r}^{(i)}(\bm x_{\cdot r}^{(i)})^\top[\widehat{\bm J}-\bm J]\},
\end{aligned}
$$
where $\widetilde{\bm J}^{(j)}$ is another point on the line segment joining $\widehat{\bm J}$ and $\bm J$. Setting $$a_i:=h^{\prime}(x_r^{(i)} (\widetilde{\bm J}^{(j)})^\top \bm x_{\cdot r}^{(i)})  x_r^{(i)}x_{\cdot r, j}^{(i)} \quad \text{and}\quad b_i:=[\bar{\bm J}^{(j)}-\bm J]^\top \bm x_{\cdot r}^{(i)}(\bm x_{\cdot r}^{(i)})^\top[\widehat{\bm J}-\bm J]~,$$ we have:
$$
\left|R_j\right| =\frac{1}{n}\left|\sum_{i=1}^n a_i b_i\right| \leq \frac{1}{n}\|\bm a\|_{\infty}\|\bm b\|_1 .
$$
A calculation shows that $\|\bm a\|_{\infty} \leq 8(k!)^3$, and
$$
\begin{aligned}
\frac{1}{n}\|\bm b\|_1 &=t_j[\widehat{\bm J}- \bm J]^T\left\{\frac{1}{n} \sum_{i=1}^n \bm x_{\cdot r}^{(i)}(\bm x_{\cdot r}^{(i)})^T\right\}[\widehat{\bm J}-\bm J] \\
&=t_j[\widehat{\bm J}_S-\bm J_S^*]^T\left\{\frac{1}{n} \sum_{i=1}^n \bm x_S^{(i)}(\bm x_S^{(i)})^T\right\}[\hat{\bm J}_S- \bm J_S]\quad(\text{since}~\widehat{\bm J}_{S^c}= \bm J_{S^c}^*=0) \\
& \leqslant D_{\max }\|\widehat{\bm J}_S - \bm J_S\|_2^2~.
\end{aligned}
$$
Combining everything, we have:
$$\|\bm R\|_\infty \le 8(k!)^3 D_{\max} \|\widehat{\bm J}_S - \bm J_S\|_2^2$$
Hence, by Lemma \ref{l2con}, we have:
$$\frac{\|\bm R\|_\infty}{\lambda} \le 50 (k!)^3 \lambda d\frac{D_{\max}}{C_{\min}^2}~.$$ Lemma \ref{l337} now follows from the hypothesis about the bound on $\lambda d$.
 
\section{Proof of Lemma \ref{lem5}}\label{aplem5}
We only prove \eqref{lem5eqn1}. The proof of \eqref{lem5eqn2} is analogous, and we skip it. Now, note that if $\bm y \in \mathbb{R}^d$ is a unit-norm minimal eigenvector of $\widehat{\bm Q}_{SS}$, then:
$$
\begin{aligned}
\Lambda_{\min }\left({\bm Q}_{S S}\right) &=\min _{\|\bm x\|_2=1} \bm x^\top {\bm Q}_{SS} \bm x \\
&=\min _{\|\bm x\|_2=1}\left\{\bm x^\top \widehat{\bm Q}_{S S} \bm x + \bm x^\top\left(\bm Q_{S S}-\widehat{\bm Q}_{S S}\right) \bm x\right\} \\
& \leq \bm y^\top \widehat{\bm Q}_{S S} \bm y+ \bm y^\top\left(\bm Q_{S S}-\widehat{\bm Q}_{S S}\right) \bm y,
\end{aligned}
$$
where $\bm y \in \mathbb{R}^d$ is a unit-norm minimal eigenvector of $\widehat{\bm Q}_{S S}$. Hence, we have

\begin{equation}\label{frstbd}
    \Lambda_{\min }(\widehat{\bm Q}_{S S}) \geq \Lambda_{\min }(\bm Q_{S S})-\vertiii{\bm Q_{S S}-\widehat{\bm Q}_{S S}}_2 \geq C_{\min }-\vertiii{\bm Q_{S S}-\widehat{\bm Q}_{S S}}_2.
\end{equation}

We now aim to bound the term $\vertiii{\bm Q_{S S}-\widehat{\bm Q}_{S S}}_2$. Towards this, note that the $(j,k)^{\mathrm{th}}$ entry of the matrix $\widehat{\bm Q} - \bm Q$ can be written as:
$$Z_{j,k} := \frac{1}{n}\sum_{i=1}^n Z_{j,k}^{(i)}$$ where $\{Z_{j,k}^{(i)}\}_{i=1}^n$ is an i.i.d. sequence of mean zero and bounded (by $8(k!)^2$) random variables. Hence, by Hoeffding's inequality, we have:
\begin{equation}\label{hoeffnews}
    \p(Z_{j,k}^2 \ge \varepsilon^2) = \p\left(\left|\sum_{i=1}^n Z_{j,k}^{(i)}\right|\ge n\varepsilon\right) \le 2\exp\left(- \frac{n\varepsilon^2}{128(k!)^4}\right)
\end{equation}
Now, note that:

$$
\vertiii{\widehat{\bm Q}_{S S}-\bm Q_{S S}}_2 \leq\left(\sum_{j\in S}\sum_{k\in S} Z_{j k}^2\right)^{1 / 2},
$$
Setting $\varepsilon^2=\delta^2 / d^2$ in \eqref{hoeffnews}, we have by a union bound,
\begin{equation}\label{lem5eqn3}
    \mathbb{P}\left(\vertiii{\widehat{\bm Q}_{S S}- \bm Q_{S S}}_2 \geq \delta\right) \leq 2 \exp \left(-\frac{n\delta^2}{128 d^2 (k!)^4}+2 \log d\right)
\end{equation}
thereby proving Lemma \ref{lem5}, in view of \eqref{frstbd}.

\section{Proof of Lemma \ref{sampinc}}\label{apsampinc}
To begin with, note that we can write: 
$$\widehat{\bm Q}_{S^c S}\widehat{\bm Q}_{S S}^{-1}=\bm T_1+ \bm T_2+ \bm T_3+ \bm T_4$$ where:
\begin{eqnarray*}
    && \bm T_1:={\bm Q}_{S^c S}[\widehat{\bm Q}_{SS}^{-1}-\bm Q_{S S}^{-1}]\\ && \bm T_2:=[\widehat{\bm Q}_{S^c S}- \bm Q_{S^c S}]\bm Q_{S S}^{-1}\\ && \bm T_3:=[\widehat{\bm Q}_{S^c S}-\bm Q_{S^c S}][\widehat{\bm Q}_{S S}^{-1}-\bm Q_{S S}^{-1}]\\ && \bm T_4:= \bm Q_{S^c S} \bm Q_{S S}^{-1}.
\end{eqnarray*}
By the population incoherence assumption \eqref{a2}, we have:
$$
\vertiii{\bm T_4}_{\infty} \leq 1-\alpha,
$$
Hence, in order to complete the proof of Lemma \ref{sampinc}, it suffices to show that $\vertiii{T_j}_\infty \le \alpha/6$ with high probability for $j \in \{1,2,3\}$. The following lemma helps us achieve this:

\begin{lem}\label{lem6p6}
For any $\delta>0$ and constants $K, K^{\prime}$, the following bounds hold:
    $$\mathbb{P}\left(\vertiii{\widehat{\bm Q}_{S^c S}- \bm Q_{S^c S}}_{\infty} \geq \delta\right)
    \leq 2 \exp \left(-K \frac{n \delta^2}{d^2 (k!)^4}+\log d +\log \left[{p-1 \choose{ k-1}}-d\right]\right)~,$$
    $$\mathbb{P}\left(\vertiii{\widehat{\bm Q}_{S S}- \bm Q_{S S}}_{\infty} \geq \delta\right) \leq 2 \exp \left(-K \frac{n \delta^2}{d^2 (k!)^4}+2 \log d\right)~,$$
    $$\mathbb{P}\left(\vertiii{\widehat{\bm Q}_{S S}^{-1} -\bm Q_{SS}^{-1}}_{\infty} \geq \delta\right)
    \leq 4 \exp \left(-K \frac{n \delta^2}{d^3 (k!)^4}+K^{\prime} \log d\right)~.$$
\end{lem} 

\begin{proof}[Proof of Lemma \ref{lem6p6}]
To begin with, note that for every $\delta >0$, we have:
$$
\begin{aligned}
    \mathbb{P}\left(\vertiii{\widehat{\bm Q}_{S^c S}- \bm Q_{S^c S}}_{\infty} \geq \delta\right) &= \mathbb{P}\left[\max_{s\in S^c} \sum_{t\in S} |Z_{s,t}|\geqslant \delta\right]\\
    &\leqslant \left[{p-1 \choose k-1}-d\right] \max_{s\in S^c}\mathbb{P}\left(\sum_{t\in S} |Z_{s,t}|\geqslant \delta\right),
\end{aligned}
$$
By a further union bound, we have
$$\mathbb{P}\left(\sum_{t\in S} |Z_{s,t}|\geqslant \delta\right)
    \le d \max_{t\in S}\mathbb{P}(Z_{s,t}|\geqslant \delta/d),$$
Combining these, we have:
$$
\mathbb{P}\left(\vertiii{\widehat{\bm Q}_{S^c S} - \bm Q_{S^c S}}_{\infty} \geq \delta\right)\leqslant \left[{p-1 \choose k-1}-d\right]d ~\max_{s\in S^c, ~t\in S}\mathbb{P}[|Z_{s, t}|\geqslant \delta/d],
$$
Setting $\epsilon = \frac{\delta}{d}$ in \eqref{hoeffnews}, we get:
$$\mathbb{P}\left(\vertiii{\widehat{\bm Q}_{S^c S} - \bm Q_{S^c S}}_{\infty} \geq \delta\right)\leqslant 2\left[{p-1 \choose k-1}-d\right]d \exp\left(-\frac{n\delta^2}{128d^2(k!)^4}\right)$$
which gives the first inequality.
The proof of the second inequality is analogous, with the factor ${p-1 \choose k-1} -d$ being replaced with $d$. To show the last inequality, note that:
$$
\begin{aligned}
    \vertiii{\widehat{\bm Q}_{S S}^{-1}- \bm Q_{SS}^{-1}}_{\infty} &= \vertiii{\bm Q_{S S}^{-1} (\bm Q_{S S}- \widehat{\bm Q}_{SS}) \widehat{\bm Q}_{S S}^{-1}}_{\infty}\\
    & \leqslant \sqrt{d}\vertiii{\bm Q_{S S}^{-1} (\bm Q_{S S}- \widehat{\bm Q}_{SS}) \widehat{\bm Q}_{S S}^{-1}}_{2}\\
    & \leqslant \sqrt{d}\vertiii{\bm Q_{S S}^{-1}}_2 \vertiii{\bm Q_{S S}- \widehat{\bm Q}_{SS}}_2\vertiii{\widehat{\bm Q}_{S S}^{-1}}_{2}\\
    & \leqslant \frac{\sqrt{d}}{C_{\min}}\vertiii{\bm Q_{S S}- \widehat{\bm Q}_{SS}}_2\vertiii{\widehat{\bm Q}_{S S}^{-1}}_{2}.
\end{aligned}
$$
From \eqref{lem5eqn3}, we have:
\begin{equation*}
    \mathbb{P}\left(\vertiii{\widehat{\bm Q}_{S S}^{-1}}_2 \geqslant \frac{2}{C_{\min}}\right) \leq 2 \exp \left(-K_1\frac{C_{\min}^2 n}{d^2 (k!)^4}+2 \log d\right)
\end{equation*}
for some constant $K_1>0$. Also, by \eqref{lem5eqn3}, we have:
\begin{equation*}
    \mathbb{P}\left(\vertiii{\widehat{\bm Q}_{S S} - \bm Q_{S S}}_2 \geqslant \frac{\delta}{\sqrt{d}}\right) \leq 2 \exp \left(-K \frac{\delta^2 n}{d^3 (k!)^4}+2 \log d\right)
\end{equation*}
for some constant $K>0$.
Combined together, we arrive at the last inequality.
\end{proof}

We are now in a position to complete the proof of Lemma \ref{sampinc}. To start with, we rewrite $\bm T_1$ as:
$$
\bm T_1= \bm Q_{S^c S} \bm Q_{S S}^{-1}\left(\bm Q_{S S}-\widehat{\bm Q}_{S S}\right)\widehat{\bm Q}_{S S}^{-1}
$$
and hence, we have:
$$
\begin{aligned}
\vertiii{\bm T_1}_{\infty} & \leq\vertiii{\bm Q_{S^c S} \bm Q_{S S}^{-1}}_\infty \vertiii{\widehat{\bm Q}_{S S}-\bm Q_{S S}}_{\infty}\vertiii{\widehat{\bm Q}_{S S}^{-1}}_{\infty} \\
& \leq (1-\alpha)\sqrt{d}\vertiii{\widehat{\bm Q}_{S S}-\bm Q_{S S}}_{\infty}\vertiii{\widehat{\bm Q}_{S S}^{-1}}_{2}
\end{aligned}
$$
On taking $\delta := C_{\min}/2$ in \eqref{lem5eqn2} of Lemma \ref{lem5}, we have:
$$\vertiii{\widehat{\bm Q}_{S S}^{-1}}_2= \left[\Lambda_{\min}\left(\widehat{\bm Q}_{S S}\right)\right]^{-1} \leq \frac{2}{C_{\min }}$$ with probability at least $1-\exp \left(-K n / d^2 (k!)^4+ B \log d\right)$ for some positive constants $K$ and $B$. Also, on taking $\delta := cd^{-1/2}$ for some positive constant $c$ in Lemma \ref{lem6p6}, we have:
$$\vertiii{\widehat{\bm Q}_{SS} - \bm Q_{SS}}_\infty \leq \frac{c}{\sqrt{d}}$$ with probability at least $1-2 \exp \left(-K n c^2 / d^3 (k!)^4+ 2\log d\right)$. Hence, $c$ can be chosen sufficiently small, to ensure that there exist constants $K'$ and $K''$ such that
$$\vertiii{\bm T_1}_\infty \le \frac{\alpha}{6}$$ with probability at least $1-3\exp \left(-K' n / d^3 (k!)^4+ K''\log d\right)$

Next, note that: 
$$
\begin{aligned}
    \vertiii{\bm T_2}_{\infty}\leqslant & \sqrt{d}\vertiii{\bm Q_{SS}^{-1}}_2\vertiii{\widehat{\bm Q}_{S^c S}- \bm Q_{S^c S}}_{\infty}\\
    \leqslant & \frac{\sqrt{d}}{C_{\min}}\vertiii{\widehat{\bm Q}_{S^c S}-\bm Q_{S^c S}}_{\infty},
\end{aligned}
$$
Applying the first bound in Lemma \ref{lem6p6} with $\delta = \frac{\alpha C_{\min}}{6\sqrt{d}}$, we have:
\begin{equation*}
    \mathbb{P}\left(\vertiii{{\bm T}_2}_{\infty} \geq \frac{\alpha}{6}\right) \leq 2 \exp \left(-K \frac{n}{d^3 (k!)^4}+\log d +\log \left[{p-1 \choose k-1}-d\right]\right)
\end{equation*}
for some constant $K>0$.

Similarly applying the first and the last bounds in Lemma \ref{lem6p6} with $\delta = \sqrt{\frac{\alpha}{6}}$, we have: 
\begin{equation*}
    \mathbb{P}\left(\vertiii{{\bm T}_3}_{\infty} \geq \frac{\alpha}{6}\right) \leq 4 \exp \left(-K_1 \frac{n }{d^3 (k!)^4}+ K_2 \log d +\log \left[{p-1 \choose k-1}-d\right]\right) 
\end{equation*}
for some constants $K_1,K_2>0$. Combining all the above bounds, we conclude the proof of Lemma \ref{sampinc} .

%\pagebreak
%\bibliographystyle{plain} % We choose the "plain" reference style
%\bibliography{ref} % Entries are in the refs.bib file

\end{document}